\documentclass[aos,MSNbibl,nameyear,dvips]{arximspdf}
\usepackage{dcolumn}
\usepackage{graphicx}

\doi{10.1214/12-AOS1063} 
\volume{40}
\issue{6}
\pubyear{2012}
\firstpage{3003}
\lastpage{3030}

\makeatletter
\newcolumntype{d}[1]{D{.}{.}{#1}}
\renewcommand{\bm}{\mathbf}
\newcommand{\bmm}{\bolds}
\newcommand{\eqref}[1]{(\ref{#1})}

\newcommand{\bA}{\mathbf{A}}
\newcommand{\bB}{\mathbf{B}}
\newcommand{\bC}{\mathbf{C}}
\newcommand{\bI}{\mathbf{I}}
\newcommand{\bS}{\mathbf{S}}
\newcommand{\bP}{\mathbf{P}}
\newcommand{\bW}{\mathbf{W}}
\newcommand{\bX}{\mathbf{X}}
\newcommand{\bY}{\mathbf{Y}}
\newcommand{\bD}{\mathbf{D}}
\newcommand{\bL}{\mathbf{L}}
\newcommand{\bM}{\mathbf{M}}

\newcommand{\bepsilon}{\bolds{\varepsilon}}
\newcommand{\be}{\mathbf{e}}
\newcommand{\bbeta}{\bolds{\beta}}
\newcommand{\bSigma}{\bolds{\Sigma}}
\newcommand{\bLambda}{\bolds{\Lambda}}
\newtheorem{theorem}{Theorem}[section]
\newtheorem{lem}{Lemma}[section]

\newproclaim{remark}{Remark}
\makeatother

\begin{document}
\begin{frontmatter}

\title{Asymptotic optimality and efficient computation of the
leave-subject-out cross-validation}
\runtitle{leave-subject-out cross-validation}

\begin{aug}
\author[A]{\fnms{Ganggang} \snm{Xu}\thanksref{t1}\ead[label=e1]{gang@stat.tamu.edu}}
\and
\author[A]{\fnms{Jianhua Z.} \snm{Huang}\thanksref{t1,t2}\corref{}\ead[label=e2]{jianhua@stat.tamu.edu}}
\runauthor{G. Xu and J.~Z. Huang}
\thankstext{t1}{Supported in part by Award Number
KUS-CI-016-04, made by King Abdullah University of Science and
Technology (KAUST).}
\thankstext{t2}{Supported
in part by NSF Grants DMS-09-07170, DMS-10-07618, DMS-12-08952 and NCI Grant CA57030.}
\affiliation{Texas A\&M University}
\address[A]{Department of Statistics\\
Texas A\&M University\\
College Station, Texas 77843-3143\\
USA\\
\printead{e1}\\
\phantom{E-mail:\ } \printead*{e2}} 
\end{aug}

\received{\smonth{4} \syear{2012}}
\revised{\smonth{9} \syear{2012}}

%
\begin{abstract}
Although the leave-subject-out cross-validation (CV)
has been widely used in practice for tuning parameter
selection for various nonparametric and semiparametric
models of longitudinal data, its theoretical property is unknown
and solving the associated optimization problem
is computationally expensive,
especially when there are multiple tuning parameters.
In this paper, by focusing on the penalized spline method,
we show that the leave-subject-out CV is optimal in the sense
that it is asymptotically equivalent to the empirical
squared error loss function minimization. An efficient
Newton-type algorithm is developed to compute the penalty
parameters that optimize the CV criterion.
Simulated and real data are used to demonstrate the effectiveness
of the leave-subject-out CV in selecting
both the penalty parameters and the working correlation matrix.
\end{abstract}

%

\end{frontmatter}

\section{Introduction}\label{sec1}
In recent years there has seen a growing interest in applying flexible
statistical models for analyzing longitudinal data or the more general
clustered data. Various semiparametric
[e.g., \citet{ZegerDiggle1994,Zhang-etal-98,LinYing2001,Wangetal2005}] and
nonparametric [e.g., \citet{FanZhang2000},
\citet{LinCarroll2000},
\citet{RiceSilverman1991},
\citeauthor{Wang98} (\citeyear{Wang98,Wang2003}),
\citet{Welshetal2002}, \citet{Zhuetal2008}]
models have been proposed
and studied in the literature. All of these flexible, semiparametric
or nonparametric methods require specification of tuning parameters,
such as
the bandwidth for the local polynomial kernel methods, the number of
knots for regression splines and the penalty parameter for
penalized splines and smoothing splines.

The ``leave-subject-out cross-validation'' (LsoCV) or more generally
called ``leave-cluster-out cross-validation,'' introduced by
\citet{RiceSilverman1991}, has been widely used as the method for
selecting tuning parameters in analyzing longitudinal data and
clustered data; see, for example,
\citet{Hooveretal1998,Huangetal2002,WuZhang06,Wangetal2008}.
The LsoCV is intuitively appealing since the within-subject
dependence is preserved by leaving out all observations from
the same subject together in the cross-validation.
In spite of its broad acceptance in practice, the use of LsoCV
still lacks a theoretical justification to date. Computationally,
the existing literature has focused on the grid search
method for finding the minimizer of the LsoCV criterion (LsoCV score)
[\citet{Chiangetal2001,Huangetal2002,Wangetal2008}], which
is rather inefficient and even prohibitive
when there are multiple tuning parameters. The goal of this paper is
twofold: First, we develop a theoretical justification of the
LsoCV by showing that the LsoCV criterion is asymptotically equivalent
to an appropriately defined loss function; second, we develop a
computationally efficient algorithm to optimize the LsoCV criterion for
selecting multiple penalty parameters for penalized splines.

We shall focus our presentation on longitudinal data, but all
discussions in this paper apply to clustered data analysis. Suppose
we have $n$ subjects and subject $i$, $i=1,\ldots,n$, has observations
$(y_{ij},\bm x_{ij})$, $j=1, \ldots, n_i$, with
$y_{ij}$ being the $j$th response and
$\bm x_{ij}$ being the corresponding vector of covariates. Denote
$\bm y_i=( y_{i1},\ldots,y_{in_i})^T$ and $\tilde{\bX}_i=(\bm x_{i1},\ldots,\bm x_{in_i})$. The marginal non- and semi-parametric
regression model [\citet{Welshetal2002,Zhuetal2008}] assumes that the
mean and covariance matrix of the responses are given by
%
\begin{equation}
\label{mean} \mu_{ij}=E(y_{ij}|\tilde{\bX}_i)=
\bm x_{ij0}\bmm\beta_0+\sum_{k=1}^mf_k(
\bm x_{ijk}), \qquad\operatorname{cov}(\bm y_i|\tilde{\bX}_i)=
\bSigma_i,
\end{equation}
where $\bmm\beta_0$ is a vector of linear regression coefficients,
$f_k$, $k=1,\ldots,m$, are unknown
smooth functions, and $\bSigma_i$'s are within-subject covariance
matrices.
Denote $\bmm\mu_i=(\mu_{i1},\ldots,\mu_{in_i})^T$.
By using a basis expansion to approximate each $f_k$, $\bmm\mu_i$ can
be approximated by $\bmm\mu_i\approx\bX_i\bmm\beta$ for some design matrix
$\bX_i$ and unknown parameter vector $\bmm\beta$, which then can be
estimated by minimizing the penalized weighted least squares
%
\begin{equation}
\label{eqpenalized} \operatorname{pl}(\bmm\beta)=\sum_{i=1}^n(
\bm y_i-\bX_i\bmm\beta)^T
\bW_i^{-1}(\bm y_i-\bX_i\bmm\beta)+
\sum_{k=1}^m\lambda_k\bmm\beta^T \bS_k\bmm\beta,
\end{equation}
where $\bW_i$'s are working correlation matrices that are
possibly misspecified, $\bS_k$ is a semi-positive definite matrix such
that $\bmm\beta^T \bS_k\bmm\beta$ serves as a roughness penalty for~$f_k$,
and $\bmm\lambda=(\lambda_1,\ldots,\lambda_m)$ is a vector
of penalty parameters.

Methods for choosing basis functions, constructing the corresponding
design matrices $\bX_i$'s and defining the roughness penalty matrices
are well established in the statistics literature. For example,
B-spline basis and basis obtained from reproducing kernel Hilbert
spaces are commonly used. Roughness penalty matrices can be formed
corresponding to the squared second-difference penalty, the squared
second derivative penalty, the thin-plate splines penalty
or using directly the reproducing kernels. We refer to the books by
\citet{GreenSilverman1994}, \citet{Gu2002} and \citet{Wood2006} for
thorough treatments of this subject.

The idea of using working correlation for longitudinal data can be
traced back to the generalized estimating equations (GEE) of
\citet{LiangZeger1986}, where it is established that the mean
function can be consistently estimated with the correct inference
even when the correlation structure is misspecified.
\citet{LiangZeger1986} further
demonstrated that using a possibly misspecified working correlation
structure $\bW$ has the potential to improve the estimation
efficiency over methods that completely ignore the
within-subject correlation. Similarly, results have been
obtained in the nonparametric setting in \citet{Welshetal2002}
and \citet{Zhuetal2008}. Commonly used working correlation structures
include compound symmetry and autoregressive models; see
\citet{Diggleetal2002} for a detailed discussion.

In the case of independent data, \citet{Li-kc1986} established
the asymptotic optimality of the generalized cross-validation
(GCV) [\citet{CravenWahba1979}] for penalty parameter selection
by showing that minimizing the GCV criterion is asymptotically
equivalent to minimizing a
suitably defined loss function. To understand the theoretical
property of LsoCV, we ask the following question in this paper:
What loss function does the LsoCV mimic or estimate and how good is
this estimation? We are able to show that the unweighted mean squared error
is the loss function that LsoCV is targeting. Specifically, we
obtain that, up to a quantity that
does not depend on the penalty parameters, the LsoCV score is
asymptotically equivalent to the mean squared error loss.
Our result provides the needed theoretical justification of
the wide use of LsoCV in practice.

In two related papers, \citet{GuMa2005}
and \citet{GuHan2008} developed modifications of the GCV for
dependent data under assumptions on the correlation structure
and established the optimality of the modified GCVs.
Although their modified GCVs work well when the correlation
structure is correctly specified up to some unknown parameters, they need
not be suitable when there is not enough prior knowledge
to make such a specification or the within-subject correlation
is too complicated to be modeled nicely with a simple structure.
The main difference between
LsoCV and these modified GCVs is that LsoCV utilizes working
correlation matrices in the estimating equations and
allows misspecification of the correlation structure.
Moreover, since the LsoCV and the asymptotic equivalent
squared error loss are not attached to
any specific correlation structure, LsoCV can be used to
select not only the penalty parameters but also the
correlation structure.

Another contribution of this paper is the development of a fast
algorithm for optimizing the LsoCV criterion. To avoid computation of
a large number of matrix inversions, we first derive
an asymptotically equivalent approximation of the LsoCV criterion
and then derive a Newton--Raphson type algorithm to optimize
this approximated criterion. The algorithm is particularly
useful when we need to select multiple penalty parameters.

The rest of the paper is organized as follows. Section~\ref{sec2} presents
the main theoretical results. Section~\ref{sec3} proposes a computationally
efficient algorithm for optimizing the LosCV criterion.
Results from some simulation studies
and a real data analysis are given in Sections~\ref{sec4} and~\ref{sec5}. All
technical proofs and computational implementations are collected
in the \hyperref[app]{Appendix} and in the supplementary materials [\citet{supp}].

%

\section{Leave-subject-out cross validation}\label{sec2}
Let $\hat{\mu}(\cdot)$ denote the estimate of the mean function
obtained by using basis expansion of unknown functions $f_k$'s
($k=1,\ldots,m$) and solving
the minimization problem (\ref{eqpenalized}) for $\bbeta$. Let $\hat
{\mu}^{[-i]}(\cdot)$
be the estimate of the mean function $\mu(\cdot)$ by the same method
but using all the data except observations from subject $i$, $1\leq
i\leq n$. The LsoCV criterion is defined as
%
\begin{equation}
\label{defLsoCV} \operatorname{LsoCV}(\bW,\bmm\lambda)= \frac{1}{n}\sum
_{i=1}^n \bigl\{\bm y_i- \hat{
\mu}^{[-i]}(\bX_i)\bigr\}^T \bigl\{\bm
y_i- \hat{\mu}^{[-i]}(\bX_i)\bigr\}.
\end{equation}
By leaving out all
observations from the same subject, the within-subject correlation is
preserved in LsoCV.
Before giving the formal justification of LsoCV, we review a
heuristic justification in Section~\ref{secheuristic}.
Section~\ref{secloss} defines the suitable loss function.
Section~\ref{secconditions} lists the regularity conditions
and Section~\ref{example} provides an example illustrating
how the regularity conditions in Section~\ref{secconditions}
can be verified using more primitive conditions.
Section~\ref{secOptimality} presents the main theoretical
result about the optimality of LsoCV.

\subsection{Heuristic justification}\label{secheuristic}
The initial heuristic justification of LsoCV by
\citet{RiceSilverman1991} is that it mimics the mean squared
prediction error (MSPE). Consider some new observations $(\bX_i, \bm
y_i^*)$, taken at the same design points as the
observed data. For a given estimator of the mean function
$\mu(\cdot)$, denoted as $\hat\mu(\cdot)$, the MSPE is defined as
\[
\operatorname{MSPE} = \frac{1}n\sum_{i=1}^n
E \bigl\|\bm y_i^* - \hat\mu (\bX_i)\bigr\|^2 =
\frac{1}{n}\operatorname{tr}(\bSigma) +\frac{1}{n}\sum_{i=1}^nE
\bigl\|\mu(\bX_i)-\hat{\mu}(\bX_i)\bigr\|^2.
\]
Using the independence between $\hat{\mu}^{[-i]}(\cdot)$ and $\bm y_i$,
we obtain that
\[
E\bigl\{\operatorname{LsoCV}(\bW,\bmm\lambda)\bigr\}=\frac{1}{n}\operatorname{tr}(\bSigma)
+\frac{1}{n}\sum_{i=1}^n E\bigl\|\mu(
\bX_i)-\hat{\mu}^{[-i]}(\bX_i)\bigr\|^2,
\]
where $\bSigma=\operatorname{diag}\{\bSigma_1,\ldots,\bSigma_n\}$. When $n$ is large,
$\hat\mu^{[-i]}(\cdot)$ should be close to $\hat\mu(\cdot)$,
the estimate that uses observations from all subjects.
Thus, we expect $E\{\operatorname{LsoCV}(\bW,\bmm\lambda)\}$ to be
close to the MSPE.

\subsection{Loss function}~\label{secloss}
We shall provide a formal justification of LsoCV by showing that the LsoCV
is asymptotically equivalent to an appropriately defined loss function.
Denote $\bm Y=(\bm y_1^T, \ldots, \bm y_n^T)^T$,
$\bX=(\bX_1^T,\ldots,\bX_n^T)^T$, and
$\bW=\operatorname{diag}\{\bW_1,\ldots, \bW_n\}$. Then, for a given choice of $\bmm
\lambda$ and
$\bW$, the minimizer of~(\ref{eqpenalized}) has a closed-form
expression
%
\begin{equation}
\label{estimator} \hat{\bmm\beta}= \Biggl(\bX^T\bW^{-1}\bX+
\sum_{k=1}^m\lambda_k
\bS_k \Biggr)^{-1}\bX^T\bW^{-1}\bm Y.
\end{equation}
The fitted mean function evaluated at the design points is given by
%
\begin{equation}
\label{eqlinear-est} \hat{\mu}(\bX|\bm Y,\bW,\bmm\lambda) = \bX\hat{\bmm\beta} =
\bA (\bW,\bmm\lambda)\bm Y,
\end{equation}
where $\bA(\bW,\bmm\lambda)$ is the hat matrix defined as
%
\begin{equation}
\label{hat} \bA(\bW,\bmm\lambda)=\bX \Biggl(\bX^T\bW^{-1}
\bX+\sum_{k=1}^m\lambda_k
\bS_k \Biggr)^{-1}\bX^T\bW^{-1}.
\end{equation}
From now on, we shall use $\bA$ for
$\bA(\bW,\bmm\lambda)$ without causing any confusion.

For a given estimator $\hat{\mu}(\cdot)$ of
$\mu(\cdot)$, define the mean squared error (MSE) loss as the true
loss function
%
\begin{equation}
\label{trueloss} L(\hat{\bmm\mu}) = \frac{1}{n}\sum
_{i=1}^n\bigl\{\hat{\mu}(\bX_i)-\mu(
\bX_i)\bigr\}^T \bigl\{\hat{\mu}(\bX_i)-\mu(
\bX_i)\bigr\}.
\end{equation}
Using \eqref{eqlinear-est}, we obtain that,
for the estimator obtained by minimizing (\ref{eqpenalized}),
the true loss function~(\ref{trueloss}) becomes
%
\begin{eqnarray}
\label{eqtrueloss} L(\bW,\bmm\lambda) &= &\frac{1}{n}(\bA\bm Y-\bmm
\mu)^T(\bA\bm Y-\bmm\mu)
\nonumber
\\[-8pt]
\\[-8pt]
\nonumber
&=& \frac{1}{n}\bmm\mu^T(\bI-\bA)^T(\bI-\bA)\bmm\mu +
\frac{1}{n}\bmm\varepsilon^T\bA^T\bA\bmm\varepsilon -
\frac{2}{n}\bmm\mu^T\bigl(\bI-\bA^T\bigr)\bA\bmm
\varepsilon,
\end{eqnarray}
where $\bmm\mu=(\mu(\bX_1)^T,\ldots,\mu(\bX_n)^T)^T$,
$\bmm\varepsilon=\bm Y-\bmm\mu$. Since $E(\bmm\varepsilon|\tilde{\bX
}_1,\ldots,\tilde{\bX}_n)=0$ and
$\operatorname{Var}(\bmm\varepsilon|\tilde{\bX}_1,\ldots,\tilde{\bX}_n)=\bSigma$, the
risk function can be derived as
%
\begin{equation}
\label{eqrisk} R(\bW,\bmm\lambda)=E\bigl\{L(\bW,\bmm\lambda)\bigr\} =
\frac{1}{n}\bmm\mu^T(\bI-\bA)^T(\bI-\bA)\bmm\mu+
\frac
{1}{n}\operatorname{tr}\bigl(\bA^T\bA\bSigma\bigr).
\end{equation}

\subsection{Regularity conditions}~\label{secconditions}
This section states some regularity conditions needed for our
theoretical results.
Noticing that unless $\bW=\bI$, $\bA$ is not symmetric. We define
a symmetric version of $\bA$ as
$\tilde{\bA}=\bW^{-1/2}\bA\bW^{1/2}$. Let $\bC_{ii}$ be the
diagonal block of $\tilde{\bA}^2$ corresponding to the $i$th
subject. With some abuse of notation (but clear from the
context), denote by $\lambda_{\mathrm{max}}(\cdot)$ and
$\lambda_{\mathrm{min}}(\cdot)$ the largest and the smallest eigenvalues
of a matrix. The regularity conditions involve the quantity
$\xi(\bSigma,\bW)=\lambda_{\mathrm{max}}(\bSigma\bW^{-1})\lambda_{\mathrm{max}}(\bW)$, which takes the minimal value $\lambda_{\mathrm{max}}(\bSigma)$ when
$\bW=\bI$ or $\bW=\bSigma$.
Let $\be_i=\bSigma_i^{-1/2}\bepsilon_i$ and $\bm u_i$ be
$n_i\times1$ vectors such that $\bm u_i^T\bm u_i=1$, $i=1,\ldots,n$.

\begin{enumerate}[\textit{Condition} 1.]
\item[\textit{Condition} 1.]
For some $K>0$,
$E\{(\bm u_i^T\be_i)^4 \} \leq K$, $ i =1, \ldots, n$.

\item[\textit{Condition} $2$.]
\begin{enumerate}[(ii)]
\item[(i)]
$\max_{1\leq i\leq n} \{\operatorname{tr}(\bA_{ii})\} = O(\operatorname{tr}(\bA)/n)=o(1)$;
\item[(ii)] $\max_{1\leq i\leq n} \{\operatorname{tr}(\bC_{ii})\} = o(1)$.
\end{enumerate}

\item[\textit{Condition} 3.]
$\xi(\bSigma,\bW)/n = o(R(\bW,\bmm
\lambda))$.

\item[\textit{Condition} 4.]
$\xi(\bSigma,\bW)\{n^{-1}\operatorname{tr}(\bA)\}^2/
\{n^{-1}\operatorname{tr}(\bA^T\bA\bSigma)\} = o(1)$.

\item[\textit{Condition} $5$.]
$\lambda_{\mathrm{max}}(\bW)\lambda_{\mathrm{max}}(\bW^{-1})O(n^{-2}\operatorname{tr}(\bA)^2)= o(1)$.
\end{enumerate}

Condition 1 is a mild moment condition that requires that each
component of the standardized residual $\be_i=\Sigma_i^{-1/2}\bmm
\varepsilon_i$
has a uniformly bounded fourth moment.
In particular, when $\bepsilon_i$'s are from the Gaussian distribution,
the condition holds with $K=3$.

Condition 2 extends the usual condition on controlling leverage,
used in theoretical analysis of linear regression models.
Note that $\{\operatorname{tr}(\bA_{ii})\}$ can be interpreted as the leverage of
subject $i$, measuring the contribution to the fit from data
of subject $i$ and $\operatorname{tr}(\bA)/n$ is the average of the leverages.
This condition says that the maximum leverage cannot be arbitrarily
larger than the average leverage or, in other words, there
should not be any dominant or extremely influential subjects.
In the special case that all subjects have the same design
matrices, the condition automatically satisfies since
$\operatorname{tr}(\bA_{ii})=\operatorname{tr}(\bA)/n$ for all $i=1,\ldots,n$. Condition~2
is likely to be violated if the $n_i$'s are very unbalanced.
For example, if $10\%$ of subjects have $20$ observations and
the rest of the subjects only have $2$ or $3$ observations each,
then $\max_{1\leq i\leq n} \{\operatorname{tr}(\bA_{ii})\}/\{n^{-1}\operatorname{tr}(\bA)\}$ can
be very large.

When $n_i$'s are bounded, any reasonable choice of $\bW$ would
generally yield a bounded value of the quantity $\xi(\bSigma,\bW)$,
and condition 3 reduces to $n R(\bW,\bmm\lambda) \to\infty$,
which simply says that the parametric rate of convergence of
risk $O(n^{-1})$ is not achievable.
This is a mild condition since we are considering
nonparametric estimation. When $n_i$'s are not bounded, condition~3's
verification should be done
on a case-by-case basis. As a special case, recent results for the
longitudinal function estimation by \citet{CaiYuan2011} indicate that
condition 3 would be satisfied in this particular setting if
$\xi(\bSigma,\bW)/n^{*}=O(1)$ and $n^{*}/n^{1/2r}\to0$ or
$\xi(\bSigma,\bW)/n^{*}=o(1)$ and $n^{*}/n^{1/2r}\to\infty$ for some
$r>1$, where $n^*=(\frac{1}{n}\sum_{i=1}^n\frac{1}{n_i})^{-1}$
is the harmonic mean of $n_1,\ldots,n_n$. This conclusion
holds for both fixed common designs and independent random designs.

Condition 4 essentially says that
$\xi(\bSigma,\bW)\{n^{-1}\operatorname{tr}(\bA)\}^2=o(R(\bW,\bmm\lambda))$.
It is straightforward to show that the left-hand
side is bounded from above by
$c(\bSigma\bW^{-1})c(\bW)\{\operatorname{tr}(\tilde{\bA})/n\}^2/\{\operatorname{tr}(\tilde{\bA
}^2)/n\}$,
where $c(\bM)=\lambda_{\mathrm{max}}(\bM)/\lambda_{\mathrm{min}}(\bM)$ is the
condition number of a matrix $\bM$. If $n_i$'s are bounded,
for choices of $\bW$ such that
$\bSigma\bW^{-1}$ and $\bW$ are not singular, to ensure
condition 4 holds it suffices to have that
$\{\operatorname{tr}(\tilde{\bA})/n\}^2/\{\operatorname{tr}(\tilde{\bA}^2)/n\}=o(1)$. For regression
splines ($\bmm\lambda=\bm0$),
this condition holds if $p/n \to0$ where $p$ is the number
of basis functions used, since $\operatorname{tr}(\tilde{\bA}^2)=\operatorname{tr}(\tilde{\bA})=p$.
For penalized splines and smoothing splines, we provide a more detailed
discussion in Section~\ref{example}.

If the working correlation matrix $\bW$ is chosen to be
well-conditioned such that its condition number
$\lambda_{\mathrm{max}}(\bW)/\lambda_{\mathrm{min}}(\bW)$ is bounded, condition~5
reduces to $\operatorname{tr}(\bA)/n \to0$, which can be verified as condition~4.

Conditions 3--5 all indicate that a bad choice of the working correlation
matrix $\bW$ may deteriorate the performance of using the LsoCV.
For example, conditions~3--5 may be violated when
$\bSigma^{-1}\bW$ or $\bW$ is nearly singular. Thus,
in practice, it is wise to avoid using a working correlation $\bW$
that is nearly singular.

We do not make the assumption that $n_i$'s are bounded.
However, $n_i$ obviously cannot grow too fast relative to
the number of subjects $n$. In particular, if $n_i$'s are too large,
$\lambda_{\mathrm{max}}(\bSigma\bW^{-1})$ can be fairly large unless
$\bW\approx\bSigma$, and $\lambda_{\mathrm{max}}(\bW)$ can be fairly large due
to increase of dimensions of the working correlation matrices for
individual subjects. Thus, conditions 3--5 implicitly impose
a limit to the growth rate of $n_i$.

\subsection{An example: Penalized splines with B-spline basis functions}
\label{example}
In this section, we provide an example where conditions 3--5 can be
discussed in a more specific manner. Consider model~(\ref{mean}) with
only one nonparametric covariate $x$ and thus there is only one penalty
parameter $\lambda$. We further assume that all eigeinvalues of
matrices $\bW$ and $\bSigma\bW^{-1}$ are bounded from below and above,
that is, there exist positive constants $c_1$ and $c_2$ such that
$c_1\leq\lambda_{\mathrm{min}}(\bW)\leq\lambda_{\mathrm{max}}(\bW)\leq c_2$ and
$c_1\leq\lambda_{\mathrm{min}}(\bSigma\bW^{-1})\leq\lambda_{\mathrm{max}}(\bSigma
\bW^{-1})\leq
c_2$. Under this assumption, it is straightforward to show that
conditions 3--5 reduce to the following conditions.
\begin{longlist}
\item[\textit{Condition} $3'$.]
$nR(\bW,\lambda)\to\infty$ as $n\to
\infty$.

\item[\textit{Condition} $4'$.]
$\{n^{-1}\operatorname{tr}(\bA)\}^2/
\{n^{-1}\operatorname{tr}(\tilde{\bA}^2)\} = o(1)$.

\item[\textit{Condition }$5'$.]
$\operatorname{tr}(\bA)/n= o(1)$.
\end{longlist}

Using Lemmas 4.1 and 4.2 from \citet{GuHan2008} and similar arguments,
we have the following three inequalities:
%
\begin{eqnarray}
\label{ine1} \operatorname{tr}\bigl\{\tilde{\bA}(c_2\lambda,\bI)\bigr\} &\leq& \operatorname{tr}
\bigl\{\tilde{\bA}(\lambda,\bW)\bigr\}\leq \operatorname{tr}\bigl\{\tilde{\bA }(c_1
\lambda,\bI)\bigr\},
\\
\label{ine2} \operatorname{tr}\bigl\{\tilde{\bA}^2(c_2\lambda,\bI)\bigr
\} &\leq& \operatorname{tr}\bigl\{\tilde{\bA}^2(\lambda,\bW)\bigr\}\leq \operatorname{tr}\bigl\{
\tilde{\bA }^2(c_1\lambda,\bI)\bigr\}
\end{eqnarray}
and
%
\begin{eqnarray}
\label{ine3} && c_1c_3^{-1}\bigl\{\bI-\tilde{
\bA}(c_2\lambda,\bI)\bigr\}
\nonumber
\\[-8pt]
\\[-8pt]
\nonumber
&&\qquad \leq\bigl\{\bI-\bA(\lambda,\bW)\bigr\}^T\bigl\{\bI-\bA(\lambda,
\bW)\bigr\} \leq c_2c_3\bigl\{\bI-\tilde{
\bA}(c_2\lambda,\bI)\bigr\},
\end{eqnarray}
where $c_3=\exp\{c_2(1+(c_1^{-1}-c_2^{-1})^2+(c_1^{-1}-c_2^{-1}))\}$.
These inequalities and the definition of the risk function
$R(\bW,\lambda)$ imply that we need only to check conditions~$3'$--$5'$ for the case that $\bW=\bI$. In particular,
\eqref{ine1}--\eqref{ine3} imply that
\begin{eqnarray*}
&& c_1c_3^{-1}\bmm\mu^T\bigl\{\bI-
\tilde{\bA}(c_2\lambda,\bI)\bigr\}^2\bmm \mu +
c_1^2\operatorname{tr}\bigl\{\tilde{\bA}^2(c_2
\lambda,\bI)\bigr\}
\\
&&\qquad \leq nR(\bW,\lambda) \leq c_2c_3\bmm\mu^T
\bigl\{\bI-\tilde{\bA}(c_1\lambda,\bI)\bigr\}^2\bmm\mu
+c_2^2\operatorname{tr}\bigl\{\tilde{\bA}^2(c_1
\lambda,\bI)\bigr\},
\end{eqnarray*}
and, therefore, to show condition $3'$, it suffices to show
%
\begin{equation}
\label{eqcond3} \bmm\mu^T\bigl\{\bI-\tilde{\bA}(\lambda,\bI)\bigr
\}^2\bmm\mu\to\infty \quad\mbox{or} \quad\operatorname{tr}\bigl\{\tilde{\bA}^2(
\lambda,\bI)\bigr\}\to\infty
\end{equation}
as $n\to\infty$.

We now use existing results from the literature to show how to
verify conditions $3'$--$5'$. Note that the notation used
in the literature of penalized splines and smoothing splines is
not always consistent. To fix notation, we denote for the rest of this
section that $\lambda^*=\lambda/N$ and
$\tilde{\bA}^*(\lambda^*)=\tilde{\bA}(\lambda,\bI)$,
where $N$ is the total number of observations from all subjects.

Let $r$ denote the order of the B-splines and consider a sequence
of knots defined on the interval $[a,b]$, $a=t_{-(r-1)}=\cdots
=t_0<t_1<\cdots<t_{K_n}<t_{K_n+1}=\cdots=t_{K_n+r}=b$. Define B-spline
basis functions recursively as
\begin{eqnarray*}
B_{j,1}(x)&=&\cases{1,&\quad $t_j\leq x<t_{j+1},$
\vspace*{2pt}
\cr
0,&\quad $\mbox{otherwise},$}
\\
B_{j,r}(x)&=&\frac{x-t_j}{t_{j+r-1}-t_j}B_{j,r-1}(x)+\frac
{t_{j+r}-x}{t_{j+r}-t_{j+1}}B_{j+1,r-1}(x)
\end{eqnarray*}
for $j=-(r-1),\ldots,K_n$. When this B-spline basis is used for basis
expansion, the $j$th row of $\bX_i$ is
$\bX_{i(j)}^T=(B_{-(r-1),r}(x_{ij}),\ldots,B_{K_n,r}(x_{ij}))$, for
$j=1,\ldots,n_i$ and $i=1,\ldots,n$. When the penalty is
the integrated squared $q$th derivative of the spline
function with $q\leq r-1$, that is, $\int(f^{(q)})^2$, the
penalty term can be written in terms of the spline coefficient
vector $\bbeta$ as
$\lambda\bbeta^T\Delta_q^TR\Delta_q\bbeta$, where $R$
is a $(K_n+r-q)\times(K_n+r-q)$ matrix with
$R_{ij}=\int_a^bB_{j,r-q}(x)B_{i,r-q}(x)\,dx$ and $\Delta_q$
is a matrix of weighted $q$th order difference operator
[\citet{Claeskensetal2009}].

We make the following assumptions: (a) $\delta=\max_{0\leq j\leq
K_n}(t_{j+1}-t_j)$ is of the order $O(K_n^{-1})$ and
$\delta/\min_{0\leq j\leq K_n}(t_{j+1}-t_j)\leq M$ for some constant
$M>0$; (b)~$\sup_{x\in[a,b]}|Q_n(x)-Q(x)|=o(K_n^{-1})$, where $Q_n$
and $Q$ are the empirical and true distribution function of all
design points $\{x_1,\ldots,x_N\}$; (c)~$K_n=o(N)$.
Define quantity\vadjust{\goodbreak} $K_q=(K_n+r-q)(\lambda^* \tilde{c}_1)^{1/(2q)}$ with
some constant $\tilde{c}_1>0$ depending on $q$ and the design
density. \citet{Claeskensetal2009} showed that, under above
assumptions, if $K_q<1$,
$\operatorname{tr}\{\tilde{\bA}^*(\lambda^*)\}$ and $\operatorname{tr}\{\tilde{\bA}^{*2}(\lambda^*)\}$
are both of the order $O(K_n)$ and
$\bmm\mu^T\{\bI-\tilde{\bA}^*(\lambda^*)\}^2\bmm\mu= O(\lambda^{*2}N
K_n^{2q}+NK_n^{-2r})$; if $K_q\geq1$, $\operatorname{tr}\{\tilde{\bA}^*(\lambda^*)\}$
and $\operatorname{tr}\{\tilde{\bA}^{*2}(\lambda^*)\}$ are of order
$O(\lambda^{*-1/(2q)})$ and
$\bmm\mu^T\{\bI-\tilde{\bA}^*(\lambda^*)\}^2\bmm\mu= O(N\lambda^{*}+
NK_n^{-2q})$.
Using these results and the results following inequalities
(\ref{ine1})--(\ref{ine3}), it is
straightforward to show that if $\lambda^*=0$ (for regression splines),
letting $K_n\to\infty$ and $K_n/n\to0$ is sufficient to guarantee
conditions $3'$--$5'$, and if $\lambda^*\neq0$ (for penalized splines),
further assuming $\lambda^*\to0$ and
$n\lambda^{*1/(2q)}\to\infty$ ensures the validity of conditions
$3'$--$5'$.

When $K_q\geq1$, the asymptotic property of the
penalized spline estimator is close to that of smoothing splines,
where the number of internal knots $K_n=N$. In fact, as discussed in
\citet{GuHan2008}, for smoothing splines, it typically holds that
$\operatorname{tr}\{\tilde{\bA}^*(\lambda^*)\}$ and $\operatorname{tr}\{\tilde{\bA}^{*2}(\lambda^*)\}$
are of order $O(\lambda^{*-1/d})$ and
$\bmm\mu^T\{\bI-\tilde{\bA}^*(\lambda^*)\}^2\bmm\mu= O(N\lambda^{*})$ for
some $d>1$ as $N\to\infty$ and $\lambda^*\to0$; see also
\citet{CravenWahba1979}, \citet{Li-kc1986} and \citet{Gu2002}.
Therefore, if one has
$\lambda^*\to0$ and $n\lambda^{*1/d}\to\infty$,
conditions $3'$--$5'$ can be verified for smoothing splines.

\subsection{Optimality of leave-subject-out CV}\label{secOptimality}
In this subsection, we provide a theoretical justification of using
the minimizer of $\operatorname{LosCV}(\bW,\bmm\lambda)$ to select
the optimal
value of the penalty parameters $\bmm\lambda$. We say that the working
correlation matrix $\bW$ is predetermined if it is determined by
observation times and/or some other covariates. One way to obtain such
$\bW$ is to use some correlation function plugged in with estimated parameters.
Naturally, it is reasonable to consider the value of $\bmm\lambda$
that minimizes the true loss function
$L(\bW,\bmm\lambda)$ as the optimal value of the penalty parameters
for a predetermined $\bW$. However, $L(\bW,\bmm\lambda)$ cannot be evaluated
using data alone since the true mean function in the definition of
$L(\bW,\bmm\lambda)$ is unknown. One idea
is to use an unbiased estimate of the risk function
$R(\bW,\bmm\lambda)$ as a proxy of $L(\bW,\bmm\lambda)$. Define
%
\begin{equation}
\label{eqU} U(\bW,\bmm\lambda)=\frac{1}{n}\bY^T(\bI-
\bA)^T(\bI-\bA)\bY +\frac{2}{n}\operatorname{tr}(\bA\bSigma).
\end{equation}
It is easy to show that
%
\begin{equation}\qquad
\label{eqUL} U(\bW,\bmm\lambda)-L(\bW,\bmm\lambda)-\frac{1}{n}
\bepsilon^T\bepsilon =\frac{2}{n}\bmm\mu^T(\bI-
\bA)^T\bepsilon -\frac{2}{n}\bigl\{\bepsilon^T\bA
\bepsilon-\operatorname{tr}(\bA\bSigma)\bigr\},
\end{equation}
which has expectation zero. Thus, if $\bSigma$ is known,
$U(\bW,\bmm\lambda) - \bepsilon^T\bepsilon/n$ is an unbiased
estimate of the risk $R(\bW,\bmm\lambda)$. Actually, the estimator
is consistent, as stated in the following theorem.

\begin{theorem}
\label{thmU} Under conditions $1$--$4$, for a predetermined $\bW$ and
a nonrandom $\bmm\lambda$, as $n\to\infty$,
\[
\label{lossrisk} L(\bW,\bmm\lambda)-R(\bW,\bmm\lambda)=o_p\bigl(R(
\bW,\bmm\lambda)\bigr)\vadjust{\goodbreak}
\]
and
\[
U(\bW,\bmm\lambda)-L(\bW,\bmm\lambda)-\frac{1}{n}\bepsilon^T
\bepsilon =o_p\bigl(L(\bW,\bmm\lambda)\bigr).
\]
\end{theorem}

This theorem shows that the function $U(\bW,\bmm\lambda)-
\bepsilon^T\bepsilon/n$, the loss function $L(\bW,\bmm\lambda)$ and
the risk function $R(\bW,\bmm\lambda)$ are asymptotically equivalent.
Thus, if $\bSigma$ is known, $U(\bW,\bmm\lambda) - \bepsilon^T\bepsilon/n$
is a consistent estimator of the risk function and, moreover,
$U(\bW,\bmm\lambda)$ can be used as a reasonable surrogate of
$L(\bW,\bmm\lambda)$ for selecting the penalty parameters,
since the $\bepsilon^T\bepsilon/n$ term does not
depend on $\bmm\lambda$. However, $U(\bW,\bmm\lambda)$ depends on
knowledge of the true covariance matrix $\bSigma$, which is
usually not available. The following result states that the LsoCV
score provides a good approximation of $U(\bW,\bmm\lambda)$, without
using the knowledge of $\bSigma$.

\begin{theorem}
\label{thmLsoCV} Under conditions $1$--$5$, for a predetermined $\bW
$ and
a nonrandom $\bmm\lambda$, as $n\to\infty$,
\[
\operatorname{ LsoCV}(\bW,\bmm\lambda) - U(\bW,\bmm\lambda)=o_p\bigl(L(\bW,\bmm
\lambda)\bigr)
\]
and, therefore,
\[
\label{eqoptscv} \operatorname{ LsoCV}(\bW,\bmm\lambda)-L(\bW,\bmm\lambda)-\frac
{1}{n}
\bepsilon^T\bepsilon =o_p\bigl(L(\bW,\bmm\lambda)\bigr).
\]
\end{theorem}

This theorem shows that minimizing $\operatorname{LsoCV}(\bW,
\bmm\lambda)$ with respect to $\bmm\lambda$ is asymptotically
equivalent to minimizing $U(\bW,\bmm\lambda)$ and is also equivalent
to minimizing the true loss function $L(\bW,\bmm\lambda)$. Unlike
$U(\bW,\bmm\lambda)$, $\operatorname{LsoCV}(\bW,\bmm\lambda)$ can
be evaluated
using the data. The theorem provides the justification of
using LsoCV, as a consistent estimator of the loss or risk function,
for selecting the penalty parameters.

\begin{remark}\label{rem1} Although the above results are presented for selection of
the penalty parameter $\bmm\lambda$ for penalized splines, the
results also hold for selection of knot numbers (or number of basis
functions) $K_n$ for regression splines when $\bmm\lambda=\bm0$
and $K_n$ is the tuning parameter to be selected.
\end{remark}

\begin{remark}\label{rem2} Since the definition of the true loss function
(\ref{trueloss}) does not depend on the working correlation
structure $\bW$, we can use this loss function to compare performances
of different choices of $\bW$, for example, compound symmetry
or autoregressive, and
then choose the best one among several candidates. Thus, the result in
Theorem~\ref{thmLsoCV} also provides a justification
for using the LsoCV to select the working correlation matrix.
This theoretical implication is also confirmed in a simulation study
in Section~\ref{seccov-struct}.
When using the LsoCV to select the working correlation matrix,
we recommend to use regression splines, that is, setting
$\bmm\lambda=\bm0$, because this choice simplifies computation
and provides more stable finite sample performance.
\end{remark}
%

\section{Efficient computation}\label{sec3}
In this section, we develop a computationally efficient
Newton--Raphson-type algorithm to minimize the LsoCV score.

\subsection{Shortcut formula}\label{sec3.1}
The definition of LsoCV would indicate that it is necessary to solve
$n$ separate minimization problems in order to find the LsoCV score.
However, a computational shortcut is available that requires solving
only one minimization problem that involves all data.
Recall that $\bA$ is the hat matrix. Let $\bA_{ii}$
denote the diagonal block of $\bA$ corresponding to the observations
of subject $i$.

\begin{lem}[(Shortcut formula)]
\label{lemLsoCVformula}
The LsoCV score satisfies
%
\begin{equation}
\label{eqLsoCV} \qquad\operatorname{ LsoCV}(\bW,\bmm\lambda)=\frac{1}{n}\sum
_{i=1}^n(\bm y_i -\hat {\bm
y}_i)^T(\bI_{ii}-\bA_{ii})^{-T}(
\bI_{ii}-\bA_{ii})^{-1}(\bm y_i -\hat{
\bm y}_i),
\end{equation}
where $\bI_{ii}$ is a $n_i\times n_i$ identity matrix, and $\hat{\bm
y}_i= \hat\mu(\bX_i)$.
\end{lem}

This result, whose proof is given in the supplementary material [\citet{supp}],
extends a similar
result for independent data [e.g., \citet{GreenSilverman1994}, page 31].
Indeed, if each subject has only one observation, then (\ref
{eqLsoCV}) reduces to $\operatorname{LsoCV}
= (1/n)\sum_{i=1}^n(y_i-\hat y_i)^2/(1-a_{ii})^2$,
which is exactly the shortcut formula for the ordinary cross-validation score.

\subsection{An approximation of leave-subject-out CV}
\label{alogrithm}
A close inspection of the short-cut formula of $\operatorname
{LsoCV}(\bW,\bmm\lambda)$
given in (\ref{eqLsoCV}) suggests that the evaluation
of $\operatorname{LsoCV}(\bW,\bmm\lambda)$ can still be computationally
expensive because of the requirement of matrix inversion and
the formulation of the hat matrix $\bA$. To further reduce the
computational cost, using Taylor's expansion
$(\bI_{ii}-\bA_{ii})^{-1} \approx\bI_{ii} + \bA_{ii}$,
we obtain the following approximation of $\rm{LsoCV}(\bW,\bmm\lambda)$:
%
\begin{equation}
\label{eqLsoCV^*} \operatorname{LsoCV}^*(\bW,\bmm\lambda)=\frac{1}{n}
\bY^T(\bI-\bA )^T(\bI-\bA)\bY+\frac{2}{n}\sum
_{i=1}^n\hat{\be}_i^T
\bA_{ii}\hat{\be}_i,
\end{equation}
where $\hat{\be}_i$ is the part of $\hat{\be}=(\bI-\bA)\bY$ corresponding
to subject $i$. The next theorem shows that this
approximation is a good one in the sense that its minimization
is asymptotically equivalent to the minimization of the true
loss function.

\begin{theorem}
\label{thmLsoCV^*} Under conditions 1--5, for a predetermined $\bW$ and
a nonrandom $\bmm\lambda$, as $n\to\infty$, we have
\[
\label{eqoptscv*} \operatorname{LsoCV}^*(\bW,\bmm\lambda)-L(\bW,\bmm\lambda)-
\frac{1}{n}\bmm \varepsilon^T\bmm\varepsilon=o_p
\bigl(L(\bW,\bmm\lambda)\bigr).
\]
\end{theorem}

This result and Theorem~\ref{thmLsoCV} together imply that
$\operatorname{LsoCV}^*(\bW,\bmm\lambda)$ and $\operatorname
{LsoCV}(\bW,\bmm\lambda)$
are asymptotically equivalent, that is,
for a predetermined\vadjust{\goodbreak} $\bW$ and a nonrandom $\bmm\lambda$,
$\operatorname{LsoCV}(\bW,\bmm\lambda)-\operatorname{LsoCV}^*(\bW
,\bmm\lambda)
=o_p(L(\bW,\bmm\lambda))$.
The proof of Theorem~\ref{thmLsoCV^*} is given in the \hyperref[app]{Appendix}.

We developed an efficient algorithm to minimizing $\operatorname
{LsoCV}^*(\bW,\bmm\lambda)$ with respect to $\bmm\lambda$ for a
pre-given $\bW$ based on the works of
\citet{GuWahba1991} and \citet{Wood2004}. The idea is to optimize
the log transform of $\bmm\lambda$ using the Newton--Raphson method.
The detailed algorithm is described in the supplementary material [\citet{supp}] and
it can be shown that, for $\operatorname{LsoCV}^*(\bW,\bmm\lambda)$,
the overall computational cost
for each Newton--Raphson iteration is $O(Np)$,
which is much smaller than the cost of directly minimizing
$\operatorname{LsoCV}(\bW,\bmm\lambda)$ ($O(Np^2)$) when the total
number of
used basis functions $p$ is large.

\section{Simulation studies}\label{sec4}
\subsection{Function estimation}~\label{secfun-est}
\begin{figure}

\includegraphics{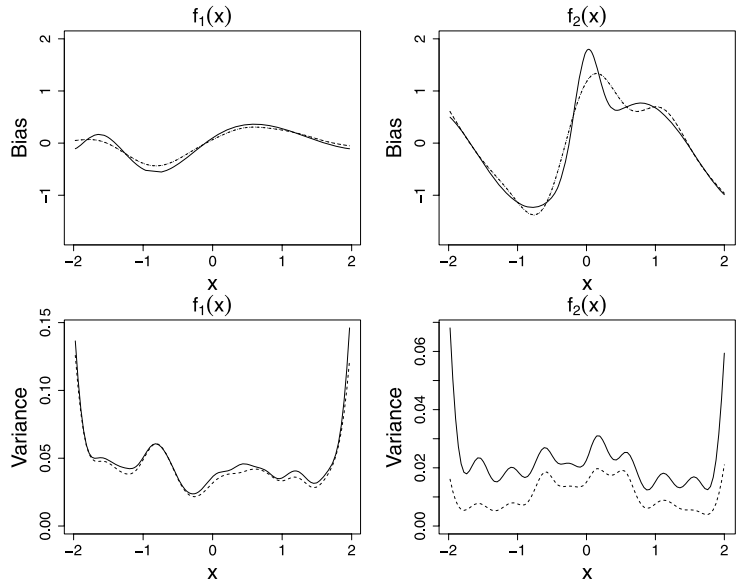}

\caption{Simulation results for function estimation based on
200 Monte Carlo runs. Functions are evaluated over $100$ equally
spaced grid points in $[-2,2]$.
Top panels: estimated functions: solid---true functions;
dashed---average of estimates using~$\bW_1$; dotted---average
of estimates using $\bW_2$ (not distinguishable with dashed).
Bottom panels: variance of estimated functions: solid---estimates
using~$\bW_1$; dashed---estimates using $\bW_2$.}
\label{fig-1}
\end{figure}

In this section, we illustrate the finite-sample performance
of $\operatorname{LsoCV}^*$ in selecting the penalty parameters.
In each simulation run, we set $n=100$ and $n_i=5$, $i=1,\ldots,
n$. A random sample is generated from the model
%
\begin{equation}
\label{sim1} y_{ij}=f_1(x_{1,i})+f_2(x_{2,ij})+
\varepsilon_{ij},\qquad  j=1,\ldots,5, i=1,\ldots, 100,
\end{equation}
where $x_1$ is a subject level covariate and $x_2$ is an
observational level covariate, both of which are drawn from
$\operatorname{Uniform}(-2,2)$. Functions used here are from \citet{Welshetal2002}:
\begin{eqnarray*}
f_1(x)&=&\sqrt{z(1-z)}\sin\biggl(2\pi\frac{1+2^{-3/5}}{1+z^{-3/5}}\biggr),
\\
f_2(x)&=&\sin(8z-4)+2\exp\bigl(-256(z-0.5)^2\bigr),
\end{eqnarray*}
where $z=(x+2)/4$. The error term $\varepsilon_{ij}$'s are generated
from a Gaussian distribution with zero mean, variance $\sigma^2$
and the compound symmetry within-subject correlation, that
is,
%
\begin{equation}
\label{corr1} \operatorname{Corr}(\varepsilon_{ij}, \varepsilon_{kl}) =
\cases{ 1, &\quad $\mbox{if $i=j=k=l$},$\vspace*{2pt}
\cr
\rho, &\quad $\mbox{if $i=k$, $j
\neq l$},$\vspace*{2pt}
\cr
0,&\quad  $\mbox{otherwise},$}
\end{equation}
$j,l=1,\ldots, 5$, $i,k=1,\ldots,100$. In this subsection, we take
$\sigma=1$ and $\rho=0.8$. A~cubic spline with 10 equally spaced
interior knots in $[-2,2]$ was used for estimating each function.
Functions were estimated
by minimizing~(\ref{eqpenalized}) with two working correlations:
the working independence (denoted as \mbox{$\bW_1=\bI$}) and the compound
symmetry with
$\rho=0.8$ (denoted as $\bW_2$). Penalty parameters were selected by
minimizing \textup{LsoCV*} defined in \eqref{eqLsoCV^*}.
The top two panels of Figure~\ref{fig-1}
show that the biases using $\bW_1$ and $\bW_2$ are almost the same, which
is consistent with the conclusion in~\citet{Zhuetal2008} that the bias of
function estimation using regression splines does not depend on the
choice of the working correlation. The bottom two panels indicate that
using the true correlation structure $\bW_2$ yields more efficient
function estimation, and the message is more clear in the estimation
of $f_2(x)$.
%

\subsection{Comparison with an existing method}\label{sec4.2}
Assuming that the structure of $\bW$ is known up to a parameter
$\gamma$ and the true covariance matrix $\bSigma$ is attained
at $\gamma=\gamma_0$, \citet{GuHan2008} proposed to
simultaneously select $\gamma$ and $\bmm\lambda$ by minimizing the
following criterion:
%
\begin{equation}
\label{hangu} \mathrm{V}^{*}(\bW,\bmm\lambda)=\log\bigl\{
\bY^T\bW^{1/2}(\bI-\tilde{\bA })^2
\bW^{1/2}\bY/N\bigr\}-\frac{1}{N}\log|\bW|+\frac{2\operatorname{tr}(\bA
)}{N-\operatorname{tr}(\bA)},\hspace*{-35pt}
\end{equation}
where $N$ is the total number of observations.
They proved that V* is asymptotically optimal
in selecting both the penalty parameter $\bmm\lambda$ and the
correlation parameter $\gamma$, provided that the within subject
correlation structure is correctly specified.
In this section, we compare the finite sample performance
of $\operatorname{LsoCV}^*$ and V* in selecting the penalty parameter
when the working correlation matrix $\bW$ is given and fixed.

\begin{figure}

\includegraphics{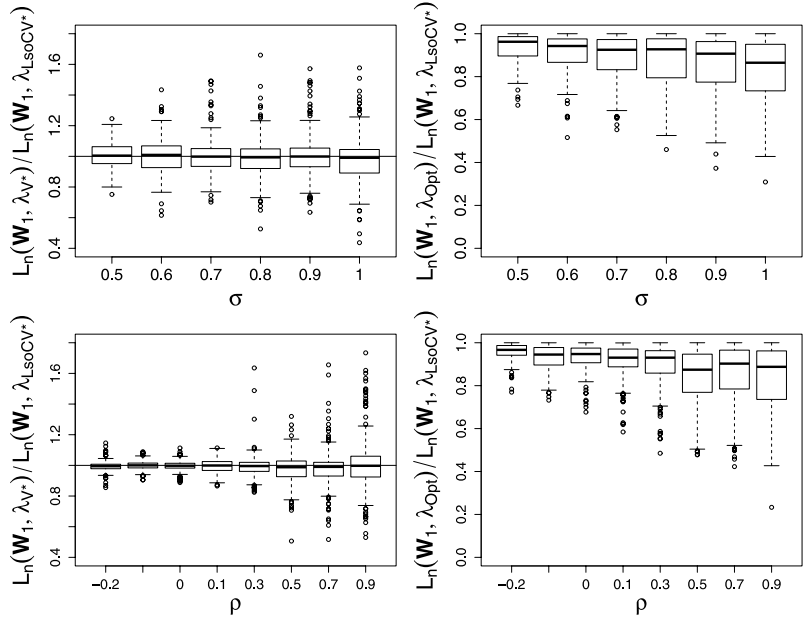}

\caption{Relative efficiency of \textup{LsoCV*} to \textup{V*} and to the
true loss when the working correlation matrix is the same as the true
correlation matrix.} \label{fig-31}
\end{figure}

We generated data using (\ref{sim1}) and (\ref{corr1})
as in the previous subsection and considered different parameters for
the correlation matrix. In particular, we fixed $\rho=0.8$ and
varied the noise standard deviation $\sigma$ from $0.5$ to $1$;
we also fixed $\sigma=1$ and varied $\rho$ from $-0.2$ to $0.9$.
A cubic spline with 10 equally spaced interior knots was used for each
unknown regression function. For each simulation run, to
compare the effectiveness of two selection criteria for a given
working correlation matrix~$\bW$, we calculated the
ratio of true losses at different choices of penalty parameters:
$L(\bW,\bmm\lambda_{\mathrm{V^{\ast}}})/L(\bW,\bmm\lambda_{\mathrm{LsoCV^{\ast}}})$
and
$L(\bW,\bmm\lambda_{\mathrm{Opt}})/L(\bW,\bmm\lambda_{\mathrm{LsoCV^{\ast}}})$,
where $\bmm\lambda_{\mathrm{V}^{*}}$ and $\bmm\lambda_{\mathrm{LsoCV^{\ast}}}$ are
penalty parameters selected by using V* and LsoCV*, respectively,
and $\bmm\lambda_{\mathrm{Opt}}$ is obtained by minimizing the true loss
function defined in~(\ref{trueloss}) assuming the mean function
$\mu(\cdot)$ is known.

\begin{figure}

\includegraphics{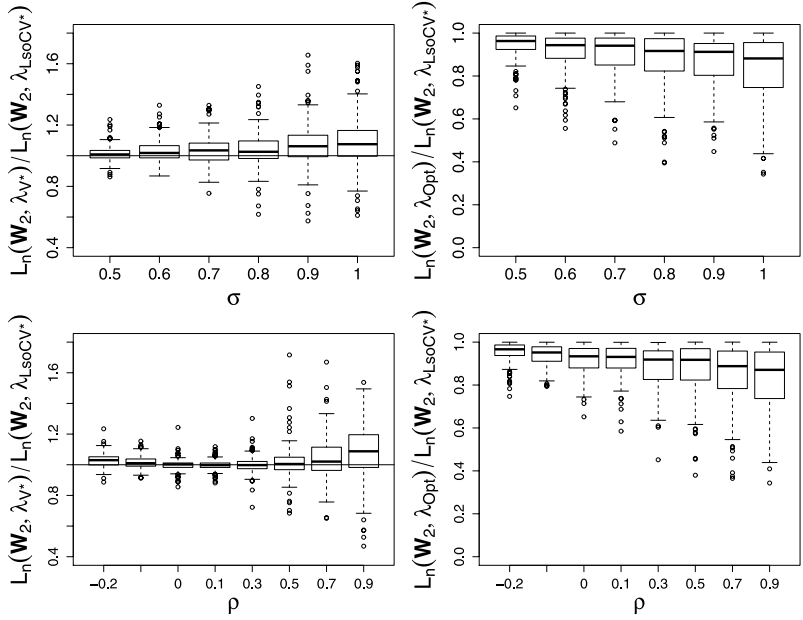}

\caption{Relative efficiency of \textup{LsoCV*} to \textup{V*} and to the
true loss when the working correlation matrix is different from the
true correlation matrix.} \label{fig-3}
\end{figure}

In the first experiment, the true correlation matrix was used as
the working correlation matrix, denoted as $\bW_1$. This is the case
that V* is
expected to work well according to \citet{GuHan2008}. Results in
Figure~\ref{fig-31} indicate that performances of LsoCV* and V*
are comparable for this case regardless of values of $\sigma$ or $\rho$.
In the second experiment, the working correlation structure was chosen
to be different from the true correlation structure. Specifically,
the working correlation matrix, denoted as $\bW_2$, is a truncated
version of (\ref{corr1})
where the correlation coefficient between $\varepsilon_{i,j_1}$ and
$\varepsilon_{i,j_2}$ is set to $\rho$ if $|j_1-j_2|=1$ and $0$ if
$|j_1-j_2|\geq2$. Results in Figure~\ref{fig-3} show that
LsoCV* becomes more effective than V* in terms of minimizing the true
loss of estimating the true mean function $\hat{\mu}(\cdot)$ as
$\sigma$ or $\rho$ increases. These results are understandable
since V* is applied to a situation that it is not designed for and its
asymptotic optimality does not hold. Moreover,
from the right two panels of
Figures~\ref{fig-31} and~\ref{fig-3}, we see that the minimum value
of LsoCV* is reasonably close to the true loss function assuming
the knowledge of the true function, as indicated by
the conclusion of Theorem~\ref{thmLsoCV^*}.

\subsection{Correlation structure selection}\label{seccov-struct}
We conducted a simulation study to evaluate
the performance of LsoCV* in selecting the working correlation matrix~$\bW$.
The data was generated using the model (\ref{sim1}) with
$\sigma=1$,
$n_i=5$ for all $i=1,\ldots,n$. In this
experiment, both $x_1$ and $x_2$ are set to be observational level
covariates drawn from $\operatorname{Uniform}(-2,2)$. Four types of within-subject
correlation structures were considered: independence (IND), compound
symmetry with correlation coefficient $\rho$ (CS), AR(1) with lag-one
correlation $\rho$ (AR), and unstructured correlation matrix with
$\rho_{12}=\rho_{23}=0.8$, $\rho_{13}=0.3$ and $0$ otherwise (UN).
Data were generated using one of these correlation structures and
then the LsoCV* was used to select the best
working correlation from the four possible candidates.
A cubic spline with $10$ equally spaced interior knots in $[-2,2]$ was used
to model each unknown function and we set the penalty parameter
vector $\bmm\lambda=\bm0$. Table~\ref{tab-1} summarizes the results
based on 200 simulation runs for each setup. We observe that
LsoCV* works well: the true correlation structure is selected in
the majority of times.

\begin{table}
\caption{Simulation results for working correlation structure
selection}\label{tab-1}
\begin{tabular*}{\textwidth}{@{\extracolsep{\fill}}lccd{3.1}d{2.1}d{2.1}d{2.1}@{}}
\hline
&&&\multicolumn{4}{c@{}}{\textbf{Selected structure}}\\[-4pt]
&&&\multicolumn{4}{c@{}}{\hrulefill}\\
$\bolds{n}$&$\bolds{\rho}$&\textbf{True structure}&\multicolumn{1}{c}{\textbf{IND}}&\multicolumn{1}{c}{\textbf{CS}}&\multicolumn{1}{c}{\textbf{AR}}&\multicolumn{1}{c@{}}{\textbf{UN}} \\
\hline
\phantom{0}50&0.3&IND&97.0&2.0&1.0&0\\
&&CS&8.5&78.0&13.5&0\\
&&AR&13.5&10.0&76.5&0\\
&&UN&1.5&1.5&21.5&75.5\\
&0.5&IND&96.5& 2.5& 1.0 & 0\\
&&CS& 3.0& 78.5& 18.5& 0\\
&&AR&4.0& 9.5& 86.5& 0\\
&&UN& 3.5& 4.0& 11.5& 81.0\\
&0.8&IND&98.5& 1.0& 0.5& 0\\
&&CS& 3.5& 74.0& 22.0& 0.5\\
&&AR&5.5& 21.0& 71.0& 2.5\\
&&UN&5.5& 1.0& 8.5& 85.0\\[3pt
]100&0.3&IND&95.0 &3.0 &2.0& 0\\
&&CS&2.0& 84.5 &13.5 &0\\
&&AR&3.5 & 8.5& 88.0& 0\\
&&UN&0& 1.0& 13.5& 85.5\\
&0.5&IND&99.5& 0.5& 0& 0\\
&&CS& 2.5& 81.0& 16.5& 0\\
&&AR&1.0 & 6.0& 93.0& 0\\
&&UN& 2.0& 0.5& 10.0& 87.5\\
&0.8&IND&99.0& 1.0& 0& 0\\
&&CS& 2.5& 73.5& 24.0 &0\\
&&AR&2.0& 20.0& 76.5& 1.5\\
&&UN&5.5& 2.0 & 9.0& 83.5\\[3pt]
150&0.3&IND& 98.5& 1.0& 0.5 & 0\\
&&CS&2.0& 85.0& 13.0& 0\\
&&AR&2.5 & 5.5& 92.0& 0\\
&&UN&0& 0& 16.5& 83.5\\
&0.5&IND&100& 0& 0 & 0\\
&&CS& 1.0& 81.5 & 17.5& 0\\
&&AR&2.5 &8.5& 89.0& 0\\
&&UN& 0.5& 0 &12.0& 87.5\\
&0.8&IND&99.5& 0.5& 0& 0\\
&&CS& 1.0& 78.0& 20.0& 1.0\\
&&AR& 0.5& 18.5& 77.5& 3.5\\
&&UN&1.0 &2.0& 6.5& 90.5\\
\hline
\end{tabular*}
\end{table}

\section{A real data example}\label{sec5}
As a subset from the Multi-center AIDS Cohort Study, the data set includes
the repeated measurements of CD4 cell counts and percentages on 283
homosexual men who became HIV-positive between 1984 and 1991.
All subjects were scheduled to take their measurements at semi-annual
visits. However, since many subjects missed some of their scheduled visits,
there are unequal numbers of repeated measurements and different
measurement times per subject.
Further details of the study can be found in \citet{Kaslowetal1987}.

\begin{figure}

\includegraphics{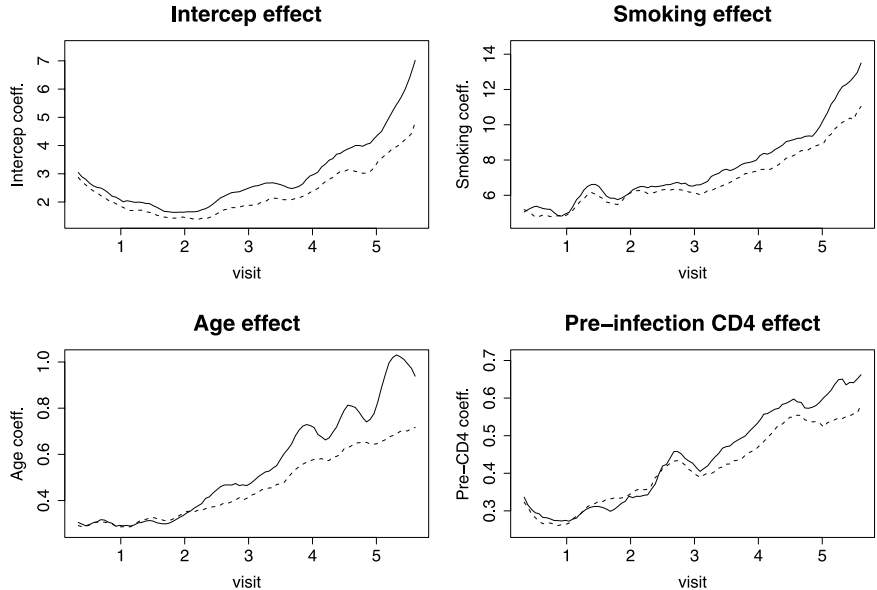}

\caption{Width of the $95\%$ pointwise bootstrap confidence
intervals based on 1000 bootstrap samples, using the working
independence $\bW_1$
(solid line) and the working correlation matrix $\bW_2$ (dashed
line).} \label{fig-4}
\end{figure}

Our goal is a statistical analysis of the trend of mean CD4
percentage depletion over time. Denote by $t_{ij}$ the time in years
of the $j$th measurement of the $i$th individual after HIV
infection, by $y_{ij}$ the $i$th individual's CD4 percentage at time
$t_{ij}$ and by $X_i^{(1)}$ the $i$th individual's smoking status
with values $1$ or $0$ for the $i$th individual ever or never smoked
cigarettes, respectively, after the HIV infection. To obtain a clear
biological interpretation, we define $X_i^{(2)}$ to be the $i$th
individual's centered age at HIV infection, which is obtained by the
$i$th individual's age at infection subtract the sample average age
at infection. Similarly, the $i$th individual's centered
pre-infection CD4 percentage, denoted by $X_i^{(3)}$, is computed by
subtracting the average pre-infection CD4 percentage of the sample
from the $i$th individual's actual pre-infection CD4 percentage.
These covariates, except the time, are time-invariant. Consider the
varying-coefficient model
%
\begin{equation}
\label{CD4} y_{ij}=\beta_0(t_{ij}) +
X_i^{(1)}\beta_1(t_{ij}) +
X_i^{(2)}\beta_2(t_{ij}) +
X_i^{(2)}\beta_2(t_{ij}),
\end{equation}
where $\beta_0(t)$ represents the trend of mean CD4 percentage
changing over time after the infection for a nonsmoker with average
pre-infection CD4 percentage and average age at HIV infection, and
$\beta_1(t)$, $\beta_2(t)$ and $\beta_3(t)$ describe the
time-varying effects on the post-infection
CD4 percentage of cigarette smoking, age at HIV infection and
pre-infection CD4 percentage, respectively. Since
the number of observations is very uneven among
subjects, we only used subjects with at least 4 observations. A cubic
spline with $k=10$ equally spaced knots was used for modeling
each function. We
first used the working independence $\bW_1=\bI$ to fit the data and then
used the residuals from this model to estimate parameters in the
correlation function
\[
\gamma(u;\alpha,\theta)=\alpha+(1-\alpha)\exp(-\theta u),
\]
where $u$ is the lag in time and $0<\alpha<1$, $\theta>0$. This
correlation function was considered previously in \citet{ZegerDiggle1994}.
The estimated parameter values are
$(\hat{\alpha},\hat{\theta})=(0.40,0.75)$. The second working
correlation matrix $\bW_2$ considered was formed using
$\gamma(u;\hat{\alpha},\hat{\theta})$. We computed that
$\operatorname{LsoCV}(\bW_1,\bm0)=881.88$ and $\operatorname
{LsoCV}(\bW_2,\bm
0)=880.33$, which implies that using $\bW_2$ is preferable.
This conclusion remains unchanged when the number of knots
varies. To visualize the gain in estimation
efficiency by using $\bW_2$ instead of $\bW_1$, we
calculated the width of the $95\%$ pointwise bootstrap confidence intervals
based on 1000 bootstrap samples, which is displayed in
Figure~\ref{fig-4}. We can observe that the bootstrap intervals using
$\bW_2$ are almost uniformly narrower than those using $\bW_1$,
indicating higher estimation efficiency.
The fitted coefficient functions (not shown to save space)
using $\bW_2$ with $\bmm\lambda$ selected by
minimizing $\operatorname{LsoCV^{\ast}}(\bW_2,\bmm\lambda)$ are similar
to those published in previous studies conducted on the same data set
[\citet{WuChiang2000,FanZhang2000,Huangetal2002}].

\begin{appendix}\label{app}
\section*{Appendix: Technical proofs}
This section is organized as follows. We first give three technical lemmas
(Lemmas~\ref{lemeigen}--\ref{lemvariance}) needed for the
proof of Theorem~\ref{thmU}. After proving Theorem~\ref{thmU}, we
give another lemma (Lemma~\ref{lemlastone}) that
facilitates proofs of Theorems~\ref{thmLsoCV} and~\ref{thmLsoCV^*}. We
prove Theorem~\ref{thmLsoCV^*} first and then proceed
to the proof of Theorem~\ref{thmLsoCV}.

Let $\lambda_{\mathrm{max}}(\bM)=\lambda_{1}(\bM)\geq\lambda_2(\bM)\geq
\cdots\geq
\lambda_{p}(\bM)=\lambda_{\mathrm{min}}(\bM)$ be eigenvalues of the $p\times p$
symmetric matrix $\bM$. We present several useful lemmas.

\begin{lem}
\label{lemeigen}
For any positive semi-definite matrices $\bM_1$ and $\bM_2$,
%
\begin{equation}
\lambda_i(\bM_1)\lambda_p(
\bM_2)\leq\lambda_i(\bM_1\bM_2)
\leq \lambda_i(\bM_1)\lambda_1(
\bM_2),\qquad i=1,\ldots, p .
\end{equation}
\end{lem}
\begin{pf}
See \citet{AndersonGupta1963} and \citet{Benasseni2002}.
\end{pf}

\begin{lem}
\label{lemtrace} For any positive semi-definite matrices $\bM_1$ and
$\bM_2$,
%
\begin{equation}
\operatorname{tr}(\bM_1\bM_2)\leq\lambda_{\mathrm{max}}(
\bM_1)\operatorname{tr}(\bM_2).
\end{equation}
\end{lem}
\begin{pf}
The proof is trivial, using the eigen
decomposition of $\bM_1$.
\end{pf}

\begin{lem}
\label{lemeigen1}
Eigenvalues of $\bA^T\bA\bSigma$ and $(\bI-\bA)^T(\bI-\bA
)\bSigma$ are bounded above by $\xi(\bSigma,\bW)=\lambda_{\mathrm{max}}(\bSigma\bW^{-1})\lambda_{\mathrm{max}}(\bW)$.
\end{lem}
\begin{pf}
Recall that $\tilde{\bA}=\bW^{-1/2}\bA\bW^{1/2}$.
For $\bA\bSigma\bA^T$, by Lemma~\ref{lemeigen},
\begin{eqnarray*}
\lambda_i\bigl(\bA^T\bA\bSigma\bigr)&=&
\lambda_i\bigl(\tilde{\bA}\bW\tilde{\bA }\bW^{-1/2}\bSigma
\bW^{-1/2}\bigr)
\\
&\leq&\lambda_i(\tilde{\bA}\bW\tilde{\bA})\lambda_{\mathrm{max}}
\bigl(\bSigma\bW^{-1}\bigr)
\\
&\leq&\lambda_i\bigl(\tilde{\bA}^2\bigr)
\lambda_{\mathrm{max}}(\bW)\lambda_{\mathrm{max}}\bigl(\bSigma\bW^{-1}
\bigr)\leq \xi(\bSigma,\bW).
\end{eqnarray*}
The last inequality follows from the fact that
$\max_{i}\{\lambda_i(\tilde{\bA}^2)\}\leq1$. Similarly, $\lambda_i((\bI-\bA)^T(\bI-\bA)\bSigma)\leq\xi(\bSigma,\bW)$ follows
from $\max_{i}\{\lambda_i((\bI-\tilde{\bA})^2)\}\leq1$.
\end{pf}

Denote
$\be=(\be^{T}_1,\ldots,\be^{T}_n)^T$, where $\be_i$'s are
independent random vectors with length $n_i$, $E(\be_i)=0$ and
$\operatorname{Var}(\be)=\bI_i$ for $i=1,\ldots,n$. For each $i$, define
$z_{ij}=(\bm u_{ij}^T\be_i)^2$ where $\bm u_{ij}^T\bm
u_{ik}=1$ if $j=k$ and 0 otherwise,
$j,k=1,\ldots,n_i$.

\begin{lem}
\label{lemvariance} If there exists a constant $K$
such that $E(z_{ij}^2)\leq K$ holds for all $j=1,\ldots,n_i$,
$i=1,\ldots,n$, then
%
\begin{equation}
\label{eqlemvariance} \operatorname{Var}\bigl(\be^{T}\bB\be\bigr)\leq2\operatorname{tr}\bigl(\bB
\bB^T\bigr) + K\sum_{i=1}^n
\bigl\{\operatorname{tr}\bigl(\bB_{ii}^*\bigr)\bigr\}^2,
\end{equation}
where $\bB$ is any $N\times N$ matrix (not necessarily symmetric),
$\bB_{ii}$ is the $i$th $(n_i\times n_i)$ diagonal block of $\bB$ and
$\bB^*_{ii}$ is an ``envelop'' matrix such that $\bB_{ii}^*\pm(\bB_{ii}+\bB_{ii}^T)/2$ are positive semi-definite.
\end{lem}

The proof of this lemma is given in the supplementary material [\citet{supp}].

\begin{pf*}{Proof of Theorem~\ref{thmU}}
In light of (\ref{eqrisk}) and (\ref{eqUL}), it suffices to show that
%
\begin{eqnarray}
\label{eq321} L(\bW,\bmm\lambda)-R(\bW,\bmm\lambda)&=&o_p\bigl(R(\bW,
\bmm\lambda)\bigr),
\\
\label{eq322} \frac{1}{n}\bmm\mu^T(\bI-\bA)^T\bmm
\varepsilon&=&o_p\bigl(R(\bW,\bmm\lambda)\bigr),
\\
\label{eq323} \frac{2}{n} \bigl\{\bepsilon^T\bA\bepsilon-\operatorname{tr}(
\bA\bSigma) \bigr\} &=&o_p\bigl(R(\bW,\bmm\lambda)\bigr)
\end{eqnarray}
because, combining (\ref{eq321})--(\ref{eq323}), we have
\[
U(\bW,\bmm\lambda)-L(\bW,\bmm\lambda)-\frac{1}{n}\bepsilon^T
\bepsilon=o_p\bigl(L(\bW,\bmm\lambda)\bigr).
\]

We first prove (\ref{eq321}). By (\ref{eqtrueloss}), we have
%
\begin{equation}
\label{eql1} \operatorname{Var}\bigl(L(\bW,\bmm\lambda)\bigr)=\frac{1}{n^2}\operatorname{Var} \bigl\{
\bepsilon^T\bA^T\bA\bepsilon-2\bmm\mu^T(\bI-
\bA)^T\bA\bepsilon \bigr\}.
\end{equation}
Define $\bB=\bSigma^{1/2}\bA^T\bA\bSigma^{1/2}$. Then
$\bepsilon^T\bA^T\bA\bepsilon=
(\bSigma^{-1/2}\bepsilon)^T\bB(\bSigma^{-1/2}\bepsilon)$. Since
$\bB$ is positive semi-definite,
by applying Lemma~\ref{lemvariance} with $\be=\bSigma^{-1/2}\bepsilon$, $\bB= \bSigma^{1/2}\bA^T\bA\bSigma^{1/2}$ and
$\bB_{ii}^*=\bB_{ii}$, we obtain
%
\begin{equation}
\label{1} \frac{1}{n^2}\operatorname{Var}\bigl(\bepsilon^T\bA^T\bA
\bepsilon\bigr)\leq\frac
{2}{n^2}\operatorname{tr}\bigl(\bB^2\bigr) +
\frac{K}{n^2}\sum_{i=1}^n\bigl\{\operatorname{tr}(
\bB_{ii})\bigr\}^2
\end{equation}
for some $K>0$ as defined in Lemma~\ref{lemvariance}. By Lemmas
\ref{lemtrace} and~\ref{lemeigen1}, under condition 3, we
have
%
\begin{eqnarray}
\label{2} \frac{2}{n^2}\operatorname{tr}\bigl(\bB^2\bigr)&\leq&
\frac{2\lambda_{\mathrm{max}}(\bA^T\bA\bSigma
)}{n^2}\operatorname{tr}\bigl(\bA^T\bA\bSigma\bigr)
\nonumber
\\[-8pt]
\\[-8pt]
\nonumber
&\leq& \frac{2\xi(\bSigma,\bW)}{n}\frac{1}{n}\operatorname{tr}\bigl(\bA^T\bA\bSigma
\bigr)=o\bigl(R^2(\bW,\bmm\lambda)\bigr).
\end{eqnarray}
Recall that $\bC_{ii}$ is the $i$th diagonal block of $\tilde{\bA}^2$.
Then, under condition 2(ii), $\operatorname{tr}(\bC_{ii})\sim o(1)$. Thus,
%
\begin{eqnarray}
\label{eql2} \operatorname{tr}(\bB_{ii})&=&\operatorname{tr}\bigl(\bL_i
\bSigma^{1/2}\bW^{-1/2}\tilde{\bA}\bW \tilde{\bA}
\bW^{-1/2}\bSigma^{1/2}\bL_i^T\bigr)\nonumber
\\
&\leq&\lambda_{\mathrm{max}}(\bW)\operatorname{tr}\bigl(\tilde{\bA}\bW^{-1/2}
\bSigma^{1/2}\bL_i^T\bL_i
\bSigma^{1/2}\bW^{-1/2}\tilde{\bA}\bigr)\nonumber
\\
&=&\lambda_{\mathrm{max}}(\bW)\operatorname{tr}\bigl(\bC_{ii}\bW_i^{-1/2}
\bSigma_i\bW_i^{-1/2}\bigr)
\\
&\leq& \lambda_{\mathrm{max}}(\bW)\lambda_{\mathrm{max}}\bigl(
\bSigma_i\bW_i^{-1}\bigr)\operatorname{tr}(\bC_{ii})\nonumber
\\
&=&o(1)\xi(\bSigma,\bW).\nonumber
\end{eqnarray}
Since
$\sum_{i=1}^n\{\operatorname{tr}(\bB_{ii})\}=\operatorname{tr}(\bB)=\operatorname{tr}(\bA^T\bA\bSigma)$, under
condition 3,
%
\begin{eqnarray}
\label{3} \frac{K}{n^2}\sum_{i=1}^n
\bigl\{\operatorname{tr}(\bB_{ii})\bigr\}^2&=&o(1)\frac{K\xi
(\bSigma,\bW)\operatorname{tr}(\bB)}{n^2}
\nonumber
\\[-8pt]
\\[-8pt]
\nonumber
&=&o(1)\frac{K\xi(\bSigma,\bW
)}{n}\frac{1}{n}\operatorname{tr}\bigl(\bA^T\bA\bSigma
\bigr)=o\bigl(R^2(\bW,\bmm\lambda)\bigr).
\end{eqnarray}
Combining (\ref{1})--(\ref{3}), we obtain
\[
\frac{1}{n^2}\operatorname{Var}\bigl(\bmm\varepsilon^T\bA^T\bA\bmm
\varepsilon\bigr)\sim o\bigl(R^2(\bW ,\bmm\lambda)\bigr).
\]
Since $\lambda_{\mathrm{max}}(\bA^T\bA\bSigma)\leq\xi(\bSigma,\bW)$ by
Lemma~\ref{lemeigen1}, under condition 3,
%
\begin{eqnarray}
\label{eql3} \frac{1}{n^2}\operatorname{Var} \bigl\{\bmm\mu^T(\bI-
\bA)^T\bA\bepsilon \bigr\} &=&\frac{1}{n^2}\bmm\mu^T(
\bI-\bA)^T\bA\bSigma\bA^T(\bI-\bA)\bmm \mu
\nonumber\\
&\leq&\frac{\lambda_{\mathrm{max}}(\bA^T\bA\bSigma)}{n}\frac{1}{n}\bmm\mu^T(\bI-
\bA)^T(\bI-\bA)\bmm\mu
\nonumber
\\[-8pt]
\\[-8pt]
\nonumber
&\leq& \frac{\xi(\bSigma,\bW)}{n}\frac{1}{n}\bmm\mu^T(\bI-\bA
)^T(\bI-\bA)\bmm\mu
\\
& =&o\bigl(R^2(\bW,\bmm\lambda)\bigr).\nonumber
\end{eqnarray}
Combining \eqref{eql1}--\eqref{eql3} and using the Cauchy--Schwarz
inequality,
we obtain $\operatorname{Var}(L(\bW,\bmm\lambda))=o(R^2(\bW,\bmm\lambda))$, which proves
(\ref{eq321}).

To show (\ref{eq322}), by Lemma (\ref{lemeigen1}) and condition
3, we have
\begin{eqnarray*}
\frac{1}{n^2}\operatorname{Var} \bigl\{\bmm\mu^T(\bI-\bA)^T\bmm
\varepsilon \bigr\} &=&\frac{1}{n^2}\bmm\mu^T(\bI-
\bA)^T\bSigma(\bI-\bA)\bmm\mu
\\
&\leq&\frac{\lambda_{\mathrm{max}}(\bSigma)}{n}\frac{1}{n}\bmm\mu^T(\bI -
\bA)^T(\bI-\bA)\bmm\mu
\\
&\leq&\frac{\xi(\bSigma,\bW)}{n}\frac{1}{n}\bmm\mu^T(\bI-\bA
)^T(\bI-\bA)\bmm\mu=o\bigl(R^2(\bW,\bmm\lambda)\bigr).
\end{eqnarray*}
The result follows from an application of the Chebyshev inequality.

To show (\ref{eq323}), applying Lemma~\ref{lemvariance} with
$\be=\bSigma^{-1/2}\bepsilon$,
$\bB=\bSigma^{1/2}\bA\bSigma^{1/2}$. For each $\bB_{ii}=\bSigma_{i}^{1/2}\bA_{ii}\bSigma_{i}^{1/2}$,
noticing that $(\bW_i^{1/2}-\alpha\bW_i^{-1/2})\tilde{\bA}_{ii}(\bW_i^{1/2}-\alpha
\bW_i^{-1/2})$ is positive semi-definite, we can define an ``envelop''
matrix as $\bB_{ii}^*=\frac{1}{2}\bSigma_i^{1/2}(\bW_i^{1/2}\times \tilde
{\bA}_{ii}\bW_{i}^{1/2}/\alpha_i+\alpha_i\bW_i^{-1/2}\tilde{\bA
}_{ii}\bW_{i}^{-1/2})\bSigma_i^{1/2}$ for any $\alpha_i>0$. Then by
Lemma~\ref{lemvariance}, we obtain
%
\begin{eqnarray}
\label{eqerror-sq1} \frac{2}{n^2}\operatorname{Var}\bigl(\bepsilon^T\bA
\bepsilon\bigr) &=&\frac{2}{n^2}\operatorname{Var}\bigl(\be^{T}\bB\be\bigr)
\nonumber
\\[-8pt]
\\[-8pt]
\nonumber
&\leq&\frac{2}{n^2}\operatorname{tr}\bigl(\bB\bB^T\bigr) + \frac{K}{n^2}
\sum_{i=1}^n\bigl\{\operatorname{tr}\bigl(
\bB_{ii}^*\bigr)\bigr\}^2,
\end{eqnarray}
where $K$ is as in Lemma~\ref{lemvariance}. By Lemma
\ref{lemtrace}, under condition 3, we have
\begin{eqnarray*}
\frac{2}{n^2}\operatorname{tr}\bigl(\bB\bB^T\bigr) & =&\frac{2}{n^2} \operatorname{tr}
\bigl(\bSigma\bA\bSigma \bA^T\bigr)\leq \frac{2\lambda_{\mathrm{max}}(\bSigma)}{n}
\frac{1}{n}\operatorname{tr}\bigl(\bA^T\bA\bSigma\bigr)
\\
&\leq&\frac{2\xi(\bSigma,\bW)}{n}\frac{1}{n}\operatorname{tr}\bigl(\bA^T\bA\bSigma
\bigr)=o\bigl(R^2(\bW,\bmm\lambda)\bigr).
\end{eqnarray*}
By using Lemma~\ref{lemeigen}
repeatedly and taking $\alpha_i=\lambda_{\mathrm{max}}(\bW_i)$, we have
\begin{eqnarray*}
\operatorname{tr}\bigl(\bB_{ii}^*\bigr)&=& \operatorname{tr}\bigl(\tilde{\bA}_{ii}
\bSigma_{i}^{1/2}\bW_i\bSigma_{i}^{1/2}
\bigr)/(2\alpha_i)+\alpha_i\operatorname{tr}\bigl(\tilde{
\bA}_{ii}\bSigma_{i}^{1/2}\bW_i^{-1}
\bSigma_{i}^{1/2}\bigr)/2
\\
&\leq&\lambda_{\mathrm{max}}\bigl(\bSigma_i\bW_i^{-1}
\bigr)\lambda_{\mathrm{max}}(\bW_i)\operatorname{tr}(\tilde{\bA}_{ii})
\\
&\leq&\xi(\bSigma,\bW)\operatorname{tr}(\tilde{\bA}_{ii}).
\end{eqnarray*}
Under conditions 2(i), 3 and 4, we have
%
\begin{equation}
\label{eqerror-sq2} \frac{K}{n^2}\sum_{i=1}^n
\bigl\{\operatorname{tr}\bigl(\bB_{ii}^*\bigr)\bigr\}^2 \leq
\frac{K}{n^2}\xi^2(\bSigma,\bW)O\bigl(n^{-2}\operatorname{tr}(
\bA)^2\bigr)=o\bigl(R^2(\bW ,\bmm\lambda)\bigr).
\end{equation}
Therefore, combining \eqref{eqerror-sq1}--\eqref{eqerror-sq2} and noticing
conditions 1--4, we have
\[
\frac{1}{n^2}\operatorname{Var}\bigl(\bepsilon^T\bA\bepsilon\bigr)\sim o
\bigl(R^2(\bW,\bmm\lambda)\bigr),
\]
which leads to \eqref{eq323}.
\end{pf*}

To prove Theorem~\ref{thmLsoCV}, it is easier to prove
Theorem~\ref{thmLsoCV^*} first. The following lemma is useful for the
proof of Theorem~\ref{thmLsoCV^*}.

\begin{lem}
\label{lemlastone} Let $\bD=\operatorname{diag}\{\bD_{11},\ldots,\bD_{nn}\}$ be
a diagonal block matrix and $\bD^*=\operatorname{diag}\{\bD_{11}^*,\ldots,\bD_{nn}^*\}$ be a positive semi-definite matrix such that $\bD^* \pm
(\bD+\bD^T)/2$ are positive semi-definite. In addition, $\bD_{ii}$'s
and\break $\bD_{ii}^*$'s meet the following conditions: \textup{(i)} $\max_{1\leq
i\leq
n}\{\operatorname{tr}(\bD_{ii}^*\bW_i)\}\sim\break
\lambda_{\mathrm{max}}(\bW)O(n^{-1}\operatorname{tr}(\bA))$; \textup{(ii)} $\max_{1\leq i\leq
n}\{\operatorname{tr}(\bD_{ii}\bW_i\bD_{ii}^T\}\sim
\lambda_{\mathrm{max}}(\bW)O(n^{-2}\operatorname{tr}(\bA)^2)$. Then, under conditions 1--5,
we have
\[
\frac{1}{n^2}\operatorname{Var} \bigl\{\bY^T(\bI-\bA)^T\bD(\bI-
\bA)\bY \bigr\} =o\bigl(R^2(\bW,\bmm\lambda)\bigr).
\]
\end{lem}

The proof is given in the supplementary material [\citet{supp}].

\begin{pf*}{Proof of Theorem~\ref{thmLsoCV^*}}
By Theorem~\ref{thmU}, it suffices to show that
\[
\operatorname{LsoCV^{\ast}}(\bW,\bmm\lambda)-U(\bW,\bmm\lambda )=o_p
\bigl(R(\bW,\bmm\lambda)\bigr),
\]
which can be obtained by showing
%
\begin{equation}
\label{L2} E \bigl\{\operatorname{LsoCV^{\ast}}(\bW,\bmm\lambda)-U(\bW,\bmm\lambda )
\bigr\}^2=o_p\bigl(R^2(\bW,\bmm\lambda)\bigr).
\end{equation}
Hence, it suffices to show that
%
\begin{eqnarray}
\label{bias} E \bigl\{\operatorname{LsoCV^{\ast}}(\bW,\bmm\lambda)-U(\bW,\bmm\lambda )
\bigr\}&=&o\bigl(R(\bW,\bmm\lambda)\bigr)\quad \mbox{and}
\\
\label{variance} \operatorname{Var}
\bigl\{\operatorname{LsoCV^{\ast}}(\bW,\bmm\lambda)-U(\bW,\bmm
\lambda ) \bigr\}&=&o\bigl(R^2(\bW,\bmm\lambda)\bigr).
\end{eqnarray}

Denote
$\bA_d=\operatorname{diag}\{\bA_{11},\ldots,\bA_{nn}\}$ and
$\tilde{\bA}_d=\operatorname{diag}\{\tilde{\bA}_{11},\ldots,\tilde{\bA}_{nn}\}
$. It
follows that $\tilde{\bA}_d=\bW^{-1/2}\bA_d\bW^{1/2}$ and
$n^{-1}\operatorname{tr}(\tilde{\bA}_d^2)=O(n^{-2}\operatorname{tr}(\bA)^2)$ by condition 2. Some
algebra yields that
\[
\operatorname{LsoCV^{\ast}}(\bW,\bmm\lambda)-U(\bW,\bmm\lambda)=\frac
{2}{n}
\bY^T(\bI-\bA)^T\bA_d(\bI-\bA)\bY-
\frac{2}{n}\operatorname{tr}(\bA \bSigma).
\]

First consider (\ref{bias}). We have that
%
\begin{eqnarray}
\label{bias1} && E \bigl\{\operatorname{LsoCV^{\ast}}(\bW,\bmm\lambda)-U(\bW,\bmm\lambda
) \bigr\}
\nonumber\\
&&\qquad =\frac{1}{n}\bmm \mu^T(\bI-\bA)^T\bigl(
\bA_d+\bA_d^T\bigr) (\bI-\bA)\bmm\mu
\\
&&\qquad\quad{} + \frac{1}{n}\operatorname{tr} \bigl\{\bA^T\bigl(\bA_d+
\bA_d^T\bigr)\bA\bSigma \bigr\}- \frac{2}{n}\operatorname{tr}
\bigl(\bA_d^T\bA_d\bSigma\bigr)-
\frac{2}{n}\operatorname{tr}\bigl(\bA_d^2\bSigma \bigr).\nonumber
\end{eqnarray}
We shall show that each term in (\ref{bias1}) is of the order
$o(R(\bW,\bmm\lambda))$.

Condition 2 says that $\max_{1\leq i\leq
n}\operatorname{tr}(\tilde{\bA}_{ii})=O(n^{-1}\operatorname{tr}(\bA))=o(1)$. Using conditions 2 and
5, we have
\begin{eqnarray*}
\operatorname{tr}\bigl(\bA_{ii}+\bA_{ii}^T\bigr)^2&=&2\operatorname{tr}
\bigl(\bA_{ii}^2+\bA_{ii}\bA_{ii}^T
\bigr)
\\
&=&2\operatorname{tr}\bigl(\tilde{\bA}_{ii}^2+\tilde{\bA}_{ii}
\bW_i\tilde{\bA }_{ii}\bW_i^{-1}
\bigr)
\\
&\leq &2\operatorname{tr}\bigl(\tilde{\bA}_{ii}^2\bigr) \bigl\{1+
\lambda_{\mathrm{max}}\bigl(\bW_i^{-1}\bigr)
\lambda_{\mathrm{max}}(\bW_i) \bigr\}
\\
&=&\lambda_{\mathrm{max}}(\bW)\lambda_{\mathrm{max}}\bigl(\bW^{-1}\bigr)O
\bigl(n^{-2}\operatorname{tr}(\bA)^2\bigr)=o(1),
\end{eqnarray*}
which implies that all eigenvalues of $(\bA_d+\bA_d^T)$ are of order
$o(1)$, and, hence,
\begin{eqnarray*}
\frac{1}{n}\bmm\mu^T(\bI-\bA)^T\bigl(
\bA_d+\bA_d^T\bigr) (\bI-\bA)\bmm\mu &=&o(1)
\frac{1}{n}\bmm\mu^T(\bI-\bA)^T(\bI-\bA)\bmm\mu =o
\bigl(R(\bW,\bmm\lambda)\bigr),
\\
\frac{1}{n}\operatorname{tr} \bigl\{\bA^T\bigl(\bA_d+
\bA_d^T\bigr)\bA\bSigma \bigr\}&=& o(1)\frac{1}{n}\operatorname{tr}
\bigl(\bA^T\bA\bSigma\bigr)=o\bigl(R(\bW,\bmm\lambda)\bigr).
\end{eqnarray*}
Under condition 4, the third term in (\ref{bias1}) can be bounded
as
%
\begin{eqnarray}
\label{14} \frac{1}{n}\operatorname{tr}\bigl(\bA_d^T
\bA_d\bSigma\bigr)&\leq&\lambda_{\mathrm{max}}\bigl(\bSigma
\bW^{-1}\bigr)\frac{1}{n}\operatorname{tr}\bigl(\tilde{\bA}_d^{1/2}
\bW^{1/2}\tilde{\bA}_d\bigr)
\nonumber\\
&\leq&\xi(\bSigma,\bW)\frac{1}{n}\operatorname{tr}\bigl(\tilde{\bA}_d^2
\bigr)
\\
&=&\xi (\bSigma,\bW)O\bigl(n^{-2}\operatorname{tr}(\bA)^2\bigr)=o\bigl(R(
\bW,\bmm\lambda)\bigr).\nonumber
\end{eqnarray}
For the last term in equation~(\ref{bias1}), observe that $(\bW_i^{1/2}-\alpha_i\bW_i^{-1/2})\bSigma(\bW_i^{1/2}-\alpha_i\bW_i^{-1/2})$ is positive semi-definite for any $\alpha_i$. Taking
$\alpha_i=\lambda_{\mathrm{max}}(\bW_i)$, we have
\begin{eqnarray*}
\frac{2}{n}\operatorname{tr}\bigl(\bA_d^2\bSigma\bigr)&=&
\frac{2}{n}\operatorname{tr}\bigl(\tilde{\bA}_d^2
\bW^{-1/2}\bSigma\bW^{1/2}\bigr)\leq\max_{1\leq i\leq n}\operatorname{tr} \bigl\{
\tilde {\bA}_{ii}^2\bigl(\bSigma_i^*+
\bSigma_i^{*T}\bigr) \bigr\}
\\
&\leq&\max_{1\leq i\leq n}\operatorname{tr} \bigl\{\tilde{\bA}_{ii}^2\bigl(
\bW_i^{1/2}\bSigma_i\bW_i^{1/2}/
\alpha_i+\alpha_i\bW_i^{-1/2}
\bSigma_i\bW_i^{-1/2}\bigr) \bigr\}
\\
& \leq&\max_{1\leq i\leq n}\bigl\{\lambda_{\mathrm{max}}\bigl(\bSigma_i
\bW_i^{-1}\bigr)\lambda_{\mathrm{max}}(\bW_i)\operatorname{tr}
\bigl(\tilde{\bA}_{ii}^2\bigr)\bigr\}
\\
&\leq&\xi(\bSigma,\bW)O\bigl(n^{-2}\operatorname{tr}(\bA)^2\bigr)=o
\bigl(R(\bW,\bmm\lambda)\bigr),
\end{eqnarray*}
where $\bSigma_i^*=\bW^{-1/2}_i\bSigma_i\bW_i^{1/2}$.
Equation~(\ref{bias1}) and thus (\ref{bias}) have been proved.\vadjust{\goodbreak}

To prove~(\ref{variance}), define $\bD=\bA_d$ and the corresponding
``envelop'' matrix $\bD^*=\operatorname{diag}\{\bD_{11}^*,\ldots,\bD_{nn}^*\}$,
where the diagonal blocks are defined as
$\bD_{ii}^*=\frac{1}{2}(\bW^{1/2}\times \tilde{\bA}_{ii}\bW_i^{1/2}/\alpha_i+\alpha_i\bW_i^{-1/2}\tilde{\bA}_{ii}\bW_i^{-1/2})$ with
$\alpha_i=\lambda_{\mathrm{max}}(\bW_i)$, then since
\begin{eqnarray*}
\operatorname{tr}\bigl(\bA_{ii}\bW_i\bA_{ii}^T
\bigr)&=&\operatorname{tr}\bigl(\tilde{\bA}_{ii}^2\bW_i\bigr)
\leq \lambda_{\mathrm{max}}(\bW_i) \bigl\{\operatorname{tr}(\bA_{ii})
\bigr\}^2\quad \mbox{and}
\\
\operatorname{tr}\bigl(\bD_{ii}^*\bW_i\bigr)&\leq&\lambda_{\mathrm{max}}(
\bW_i)\operatorname{tr}(\bA_{ii}),
\end{eqnarray*}
we have that $\max_{1\leq i\leq n}\operatorname{tr}(\bA_{ii}\bW_i\bA_{ii}^T)=\lambda_{\mathrm{max}}(\bW)O(n^{-2}\operatorname{tr}(\bA)^2)$ and that $\max_{1\leq i\leq n}\operatorname{tr}(\bD_{ii}^*\bW_i)=\lambda_{\mathrm{max}}(\bW
)O(n^{-1}\operatorname{tr}(\bA))$
by condition 2. Under conditions~3--4, (\ref{variance}) follows from
Lemma~\ref{lemlastone}.
\end{pf*}

\begin{pf*}{Proof of Theorem~\ref{thmLsoCV}}
By Theorem~\ref{thmLsoCV^*}, it suffices to show
\[
\operatorname{LsoCV}(\bW,\bmm\lambda)- \operatorname{LsoCV^{\ast}}(\bW,\bmm
\lambda)=o_p\bigl(L(\bW,\bmm\lambda)\bigr),
\]
which can be proved by showing that
\[
E \bigl\{\operatorname{LsoCV}(\bW,\bmm\lambda)-\operatorname
{LsoCV^{\ast}}(\bW,\bmm
\lambda) \bigr\}^2=o_p\bigl(R^2(\bW,\bmm
\lambda)\bigr).
\]
It suffices to show
%
\begin{eqnarray}
\label{bias3} E \bigl\{\operatorname{LsoCV}(\bW,\bmm\lambda)-\operatorname
{LsoCV^{\ast}}(\bW,\bmm\lambda) \bigr\}&=&o\bigl(R(\bW,\bmm\lambda)\bigr) \quad\mbox{and}
\\
\label{variance2} \operatorname{Var} \bigl\{\operatorname{LsoCV}(\bW,\bmm\lambda)-\operatorname
{LsoCV^{\ast}}(\bW,\bmm\lambda) \bigr\}&=&o\bigl(R^2(\bW,\bmm\lambda)\bigr).
\end{eqnarray}

For each $i=1,\ldots,n$, consider the eigen-decomposition
$\tilde{\bA}_{ii}=\bP_i\bLambda_i\bP_i^T$, where $\bP_i$ is a
$n_i\times n_i$ orthogonal matrix and
$\bLambda_i=\operatorname{diag}\{\lambda_{i1},\ldots,\lambda_{in_i}\}$,
$\lambda_{ij}\geq0$. Using this decomposition, we have
\[
(\bI_{ii}-\bA_{ii})^{-1}=\bW_i^{1/2}
\bP_i\bLambda_i^*\bP_i^T
\bW^{-1/2},
\]
where $\bLambda_i^{\ast}$ is a diagonal matrix with diagonal
elements $(1-\lambda_{ij})^{-1}$, $j=1,\ldots,n_i$. Since under
condition 2 $\max_{1\leq j\leq n_i}\{\lambda_{ij}\}\sim o(1)$, we
have $(1-\lambda_{ij})^{-1}=\sum_{k=0}^{\infty}\lambda_{ij}^k$,
which leads to
\[
(\bI_{ii}-\tilde{\bA}_{ii})^{-1}=\sum
_{k=0}^{\infty}\bP_i\bLambda_i^k
\bP_i^T=\sum_{k=0}^{\infty}
\tilde{\bA}_{ii}^k.
\]
Define
$\tilde{\bD}^{(m)}=\operatorname{diag}\{\tilde{\bD}^{(m)}_{11},\ldots,\tilde{\bD
}^{(m)}_{nn}\}$, where
$\tilde{\bD}^{(m)}_{ii}=\sum_{k=m}^{\infty}\tilde{\bA}_{ii}^k$
$i=1,\ldots,n$, $m=1,2,\ldots.$ It follows that, for each $i$,
\[
\operatorname{tr}\bigl(\tilde{\bD}_{ii}^{(m)}\bigr)= \sum
_{k=m}^{\infty}\operatorname{tr}\bigl(\tilde{\bA}_{ii}^k
\bigr)\leq \sum_{k=m}^{\infty} \bigl\{\operatorname{tr}(\tilde{
\bA}_{ii}) \bigr\}^k=\frac
{ \{\operatorname{tr}(\tilde{\bA}_{ii}) \}^m}{1-\operatorname{tr}(\tilde{\bA}_{ii})}.
\]
Since condition 2(i) gives $\max_{1\leq i\leq n}\operatorname{tr}(\bA_{ii})\sim
O(n^{-1}\operatorname{tr}(\bA))$, we obtain that
%
\begin{equation}
\label{11} \max_{1\leq i\leq
n}\operatorname{tr}\bigl(\tilde{\bD}_{ii}^{(m)}
\bigr)=O\bigl(n^{-m}\operatorname{tr}(\bA)^m\bigr),\qquad m=1,2,\ldots.
\end{equation}
Some algebra yields
\[
\operatorname{LsoCV}(\bW,\bmm\lambda)- \operatorname{LsoCV^{\ast}}(\bW,\bmm\lambda)=
\frac{1}{n}\bY^T(\bI-\bA )^T\bigl(
\bD^{(1)}+\bD^{(2)}\bigr)^{1/2}(\bI-\bA)\bY,
\]
where
$\bD^{(1)}=\bW^{-1/2}\tilde{\bD}^{(1)}\bW\tilde{\bD}^{(1)}\bW^{-1/2}$
and $\bD^{(2)}=\bW^{1/2}\tilde{\bD}^{(2)}\bW^{-1/2}$.

To show~(\ref{bias3}), note that
%
\begin{eqnarray}
\label{bias2} &&E \bigl\{\operatorname{LsoCV}(\bW,\bmm\lambda)-\operatorname
{LsoCV^{\ast}}(\bW,\bmm\lambda) \bigr\}\nonumber
\\
&&\qquad =\frac{1}{n}\bmm \mu^T(\bI-\bA)^T
\bD^{(1)}(\bI-\bA)\bmm\mu+ \frac{1}{n}\operatorname{tr} \bigl\{(\bI-
\bA)^T\bD^{(1)}(\bI-\bA)\Sigma \bigr\}
\\
&&\qquad\quad{} + \frac{1}{n}\bmm\mu^T(\bI-\bA)^T
\bD^{(2)}(\bI-\bA)\bmm\mu+ \frac{1}{n}\operatorname{tr} \bigl\{(\bI-
\bA)^T\bD^{(2)}(\bI-\bA)\Sigma \bigr\}.\nonumber
\end{eqnarray}
Using Lemmas~\ref{lemeigen} and~\ref{lemtrace} repeatedly and
condition 5, we have
\[
\lambda_{\mathrm{max}}\bigl(\bD^{(1)}\bigr)\leq\lambda_{\mathrm{max}}(
\bW)\lambda_{\mathrm{max}}\bigl(\bW^{-1}\bigr)O\bigl(n^{-2}\operatorname{tr}(
\bA)^2\bigr)=o(1).
\]
Thus, the first terms (\ref{bias2}) can be bounded as
\[
\frac{1}{n}\bmm\mu^T(\bI-\bA)^T\bD^{(1)}(
\bI-\bA)\bmm \mu=o(1)\frac{1}{n}\bmm\mu^T(\bI-\bA)^T(
\bI-\bA)\bmm \mu=o\bigl(R(\bW,\bmm\lambda)\bigr).
\]
Using Lemma~\ref{lemeigen1}, under condition 4 and (\ref{11}), the
second term of \eqref{bias2}
can be bounded as
\begin{eqnarray*}
\frac{1}{n}\operatorname{tr} \bigl\{(\bI-\bA)^T\bD^{(1)}(\bI-\bA)
\Sigma \bigr\} &\leq&\xi(\bSigma,\bW)\frac{1}{n}\operatorname{tr}\bigl(\tilde{
\bD}^{(1)2}\bigr)
\\
&=&\xi(\bSigma,\bW)O\bigl(n^{-2}\operatorname{tr}(\bA)^2\bigr)=o\bigl(R(
\bW,\bmm\lambda)\bigr).
\end{eqnarray*}
Now consider the third term in (\ref{bias2}). Under condition 5 and
(\ref{11}),
%
\begin{eqnarray}
\label{15} \operatorname{tr} \bigl\{\bigl(\bD_{ii}^{(2)}+
\bD_{ii}^{(2)T}\bigr)^2 \bigr\}&=&2\operatorname{tr}\bigl(\tilde{\bD
}_{ii}^{(2)2}\bigr)+2\operatorname{tr}\bigl(\bD_{ii}^{(2)}
\bD_{ii}^{(2)T}\bigr)
\nonumber\\
&=&2\operatorname{tr}\bigl(\tilde{\bD }_{ii}^{(2)2}\bigr)+2\operatorname{tr}\bigl(\tilde{
\bD}_{ii}^{(2)}\bW_i^{-1}\tilde{\bD
}_{ii}^{(2)}\bW_i\bigr)
\nonumber
\\[-8pt]
\\[-8pt]
\nonumber
&\leq&2\operatorname{tr}\bigl(\tilde{\bD}_{ii}^{(2)2}\bigr)+2
\lambda_{\mathrm{max}}\bigl(\bW_i^{-1}\bigr)
\lambda_{\mathrm{max}}(\bW_i)\operatorname{tr}\bigl(\tilde{\bD}_{ii}^{(2)2}
\bigr)
\\
&=&o\bigl(n^{-2}\operatorname{tr}(\bA)^2\bigr),\nonumber
\end{eqnarray}
which implies that all eigenvalues of
$\bD_{ii}^{(2)}+\bD_{ii}^{(2)T}$ are of the order
$O(n^{-1}\operatorname{tr}(\bA))$, and thus $o(1)$. Then, under conditions 1--5, we
have
\begin{eqnarray*}
\frac{1}{n}\bmm\mu^T(\bI-\bA)^T\bD^{(2)}(
\bI-\bA)\bmm \mu&=&\frac{1}{2n}\bmm \mu^T(\bI-\bA)^T
\bigl(\bD^{(2)}+\bD^{(2)T}\bigr) (\bI-\bA)\bmm\mu
\\
&=&o(1)\frac{1}{n}\bmm \mu^T(\bI-\bA)^T(\bI-\bA)\bmm
\mu=o\bigl(R(\bW,\bmm\lambda)\bigr).
\end{eqnarray*}
To study the the fourth term in (\ref{bias2}), we have
%
\begin{eqnarray}
\label{12} &&\frac{1}{n}\operatorname{tr} \bigl\{(\bI-\bA)^T
\bD^{(2)}(\bI-\bA)\bSigma \bigr\}\nonumber\\[-3pt]
&&\qquad=\frac{1}{n}\sum
_{i=1}^n\operatorname{tr} \bigl\{(\bI_{ii}-
\bA_{ii})^T\bD_{ii}^{(2)}(
\bI_{ii}-\bA_{ii})\bSigma_i \bigr\}
\\[-3pt]
&&\qquad\quad{} -\frac{1}{n}\sum_{i=1}^n\operatorname{tr}\bigl(
\bA_{ii}^T\bD_{ii}^{(2)}
\bA_{ii}\bSigma_i\bigr)+\frac{1}{n}\operatorname{tr}\bigl(
\bA^T\bD^{(2)}\bA\bSigma\bigr).\nonumber
\end{eqnarray}
To bound the first term in (\ref{12}), we note that
\begin{eqnarray*}
&&\operatorname{tr} \bigl\{(\bI_{ii}-\bA_{ii})^T
\bD_{ii}^{(2)}(\bI_{ii}-\bA_{ii})
\bSigma_i \bigr\}
\\[-3pt]
&&\qquad=\frac{1}{2}\operatorname{tr} \bigl\{(\bI_{ii}-\bA_{ii})^T
\bigl(\bW_i^{1/2}\tilde {\bD}_{ii}^{(2)}
\bW_i^{-1/2}+\bW_i^{-1/2}\tilde{
\bD}^{(2)}_{ii}\bW_i^{1/2}\bigr) (
\bI_{ii}-\bA_{ii})\bSigma_i \bigr\},
\end{eqnarray*}
which is bounded by
\begin{eqnarray*}
&&\frac{1}{2}\operatorname{tr} \bigl\{(\bI_{ii}-\bA_{ii})^T
\bigl(\bW_{i}^{1/2}\tilde {\bD}_{ii}^{(2)}
\bW_{i}^{1/2}/\alpha_i+\alpha_i
\bW_{i}^{-1/2}\tilde{\bD}_{ii}^{(2)}
\bW_{i}^{-1/2}\bigr) (\bI_{ii}-\bA_{ii})
\bSigma_{i} \bigr\}
\\[-3pt]
&&\qquad\leq\frac{1}{2}\xi(\bSigma_i,\bW_i)\operatorname{tr}\bigl(
\tilde{\bD }_{ii}^{(2)}\bigr)+\frac{\alpha_i}{2}\operatorname{tr} \bigl\{
\bigl(\tilde{\bD }_{ii}^{(2)}-2\tilde{\bD}_{ii}^{(3)}
\bigr)\bW_{i}^{-1/2}\bSigma_{i}
\bW_{i}^{-1/2} \bigr\}
\\[-3pt]
&&\qquad\quad{} +\frac{\alpha_i}{2}\operatorname{tr}\bigl\{\tilde{\bD}_{ii}^{(2)}\tilde {
\bA}_{ii}\bW_i^{-1/2}\bSigma_i
\bW_i^{-1/2}\tilde{\bA}_{ii}\bigr\}
\\[-3pt]
&&\qquad\leq\frac{1}{2}\xi(\bSigma,\bW)\bigl\{2+\lambda_{\mathrm{max}}\bigl(
\tilde{\bA }_{ii}^2\bigr)\bigr\}\operatorname{tr}\bigl(\tilde{
\bD}_{ii}^{(2)}\bigr)
\\[-3pt]
&&\qquad=o\bigl(R(\bW,\bmm\lambda)\bigr),
\end{eqnarray*}
where we take $\alpha_i=\lambda_{\mathrm{max}}(\bW_i)$. The last equation
follows from (\ref{11}) and condition~4. Similarly, we can show that
the second part of (\ref{12}) is $o(R(\bW,\bmm\lambda))$.

Consider the third part of (\ref{12}),
$\frac{1}{n}\operatorname{tr}(\bA^T\bD^{(2)}\bA\bSigma)=o(1)\frac{1}{n}\operatorname{tr}(\bA^T\bA\bSigma)= o(R(\bW,\bmm
\lambda))$ since all eigenvalues of
$\bD_{ii}^{(2)}+\bD_{ii}^{(2)T}$ are of the order $o(1)$ as is shown
in (\ref{15}). Hence, (\ref{12}) gives
\[
\frac{1}{n}\operatorname{tr}\bigl\{(\bI-\bA)^T\bD^{(2)}(\bI-\bA)
\Sigma\bigr\}=o\bigl(R(\bW ,\bmm\lambda)\bigr).
\]
Therefore, (\ref{bias3}) has been proved.

Next, we proceed to prove (\ref{variance2}). Define envelop matrices
$\bD^{(1)*}=\bD^{(1)}$ and $\bD^{(2)^*}=\operatorname{diag}\{\bD^{(2)^*}_{11},\ldots,\bD^{(2)^*}_{nn}\}$, where
$\bD^{(2)^*}_{ii}=\frac{1}{2}(\bW_i^{1/2}\tilde{\bD}_{ii}^{(2)}\bW_i^{1/2}/\alpha_i+\alpha_i\bW_i^{-1/2}\times \tilde{\bD}_{ii}^{(2)}\bW_i^{-1/2})$ with $\alpha_i=\lambda_{\mathrm{max}}(\bW_i)$. It is easy to
check that $\bD^{(1)*}$ and $\bD^{(2)*}$ are valid envelops of $\bD^{(1)}$ and $\bD^{(2)}$, respectively. Since under condition
5, we have
\begin{eqnarray*}
&&\operatorname{tr}\bigl(\bD_{ii}^{(1)}\bW_i
\bD^{(1)T}_{ii}\bigr)\\[-2pt]
&&\qquad\leq\lambda_{\mathrm{max}}(\bW )
\lambda_{\mathrm{max}}(\bW)\lambda_{\mathrm{max}}\bigl(\bW^{-1}\bigr)
\lambda_{\mathrm{max}}^2\bigl(\tilde {\bD}_{ii}^{(1)}
\bigr)\operatorname{tr}\bigl(\tilde{\bD}_{ii}^{(1)2}\bigr)
\\[-2pt]
&&\qquad = \bigl\{ \lambda_{\mathrm{max}}(\bW)\lambda_{\mathrm{max}}\bigl(
\bW^{-1}\bigr)O\bigl(n^{-2}\operatorname{tr}(\bA)^2\bigr) \bigr\}
\lambda_{\mathrm{max}}(\bW)O\bigl(n^{-2}\operatorname{tr}(\bA)^2\bigr)
\\[-2pt]
&&\qquad =\lambda_{\mathrm{max}}(\bW)O\bigl(n^{-2}\operatorname{tr}(\bA)^2\bigr),
\\[-2pt]
&&\operatorname{tr}\bigl(\bD_{ii}^{(1)*}\bW_i\bigr)\leq
\lambda_{\mathrm{max}}(\bW_i)\operatorname{tr}\bigl(\tilde{\bD }^{(1)2}_{ii}
\bigr)=\lambda_{\mathrm{max}}(\bW)O\bigl(n^{-2}\operatorname{tr}(\bA)^2
\bigr)
\end{eqnarray*}
and
\begin{eqnarray*}
\operatorname{tr}\bigl(\bD_{ii}^{(2)}\bW_i\bD^{(2)T}_{ii}
\bigr)&\leq& \lambda_{\mathrm{max}}(\bW_i)\operatorname{tr}\bigl(\tilde{
\bD}_{ii}^{(2)2}\bigr)\\
&=&\lambda_{\mathrm{max}}(\bW )O
\bigl(n^{-4}\operatorname{tr}(\bA)^4\bigr)
\\
&=&\lambda_{\mathrm{max}}(\bW)o\bigl(n^{-2}\operatorname{tr}(\bA)^2\bigr),
\\
\operatorname{tr}\bigl(\bD_{ii}^{(2)*}\bW_i\bigr)&\leq&
\lambda_{\mathrm{max}}(\bW_i)\operatorname{tr}\bigl(\tilde{\bD }_{ii}^{(2)}
\bigr)\\
&=&\lambda_{\mathrm{max}}(\bW)O\bigl(n^{-2}\operatorname{tr}(\bA)^2
\bigr).
\end{eqnarray*}
By applying Lemma~\ref{lemlastone}, we have
\[
\frac{1}{n^2}\operatorname{Var} \bigl\{\bY^T(\bI-\bA)^T
\bD^{(m)}(\bI-\bA)\bY \bigr\}=o_p\bigl(R^2(\bW,
\bmm\lambda)\bigr),\qquad m=1,2,
\]
and (\ref{variance2}) follows by the Cauchy--Schwarz inequality.
\end{pf*}
\end{appendix}

\begin{supplement}
\stitle{Efficient algorithm and additional proofs}
\slink[doi]{10.1214/12-AOS1063SUPP} 
\sdatatype{.pdf}
\sfilename{aos1063\_supp.pdf}
\sdescription{In the Supplementary Material, we give a detailed description of the algorithm proposed in Section~\ref{alogrithm}.
In addition, proofs of some technical lemmas are also included.}
\end{supplement}

%

\printaddresses

\end{document}